 \documentclass[11pt]{amsart}
 \usepackage[dvips]{graphicx}
 
\usepackage{amscd,amsmath,amsopn,amssymb,amsthm,multicol}
\usepackage[color,matrix, all, 2cell]{xy}

\usepackage{setspace}
\usepackage{upgreek}
\usepackage{textgreek}
\usepackage{enumerate}
\usepackage{color}
\usepackage{lscape}
\usepackage{tikz}
\usepackage{multirow}
\usepackage{cancel}
\usepackage{soul}
\usepackage{comment}
\usepackage{wasysym}
\usepackage{mathrsfs}
\usepackage{young}
\usepackage{mathtools}
\usepackage{bbm}
\usepackage{marginnote}
\usepackage{diagbox}
\usepackage{tikz-cd}
\usepackage{mathrsfs}
\usepackage{bbold}
\usepackage{marginnote}

\advance\oddsidemargin by -1.0cm
\advance\evensidemargin by -1.0cm
\advance\topmargin by -1.0cm

\numberwithin{equation}{section}
 \newcommand{\Ja}{\mathbb{J}}
 
 \newcommand{\R}{\mathbb{R}}
\newcommand{\Z}{\mathbb{Z}}
\newcommand{\C}{\mathbb{C}}

\DeclareMathOperator{\CP}{\C\mathsf{P}}

\DeclareMathOperator{\Hol}{\mathsf{Hol}}

\DeclareMathOperator{\Imm}{\mathsf{Im}}

\DeclareMathOperator{\Stab}{\mathsf{Stab}}
\DeclareMathAlphabet{\mathscrbf}{OMS}{mdugm}{b}{n}

\DeclareMathOperator{\SO}{\mathsf{SO}}

 \DeclareMathOperator{\SU}{\mathsf{SU}}

\DeclareMathOperator{\U}{\mathsf{U}}

\DeclareMathOperator{\G}{\mathsf{G}}
\DeclareMathOperator{\Lie}{\mathsf{Lie}}
\DeclareMathAlphabet{\mathpzc}{OT1}{pzc}{m}{it}

\DeclareMathOperator{\BLY}{\mathsf{BLY}}

\DeclareMathOperator{\Id}{\mathsf{Id}}

\newcommand{\fr}{\mathfrak}
\newcommand{\al}{\alpha}

\newcommand{\mc}{\mathcal}

\newcommand{\ep}{\varepsilon}

\newcommand{\om}{\omega}
\newcommand{\Om}{\Omega}

\usepackage{amsmath}
\DeclareFontFamily{U}{mathx}{}
\DeclareFontShape{U}{mathx}{m}{n}{<-> mathx10}{}
\DeclareSymbolFont{mathx}{U}{mathx}{m}{n}
\DeclareMathAccent{\widehat}{0}{mathx}{"70}
\DeclareMathAccent{\widecheck}{0}{mathx}{"71}

\DeclareMathAlphabet{\mathscrbf}{OMS}{mdugm}{b}{n}

\newcommand{\scF}{\mathcal{F}}

\newcommand{\scH}{\mathscr{H}}

\newcommand{\scV}{\mathscr{V}}

\newcommand{\Gg}{\ensuremath{\mathsf{G}}}

\newcommand{\Xg}{\ensuremath{\mathsf{X}}}
\newcommand{\Ss}{\ensuremath{\mathsf{S}}}

\DeclareMathOperator{\Tg}{\mathsf{T}}

\DeclareMathOperator{\Ric}{\mathsf{Ric}}
\DeclareMathOperator{\Scal}{\mathsf{Scal}}

\DeclareMathOperator{\Ker}{\mathsf{Ker}}

\DeclareMathOperator{\dd}{d}

\newtheorem{theorem}{Theorem}[section]
\newtheorem{lem}[theorem]{Lemma}
\newtheorem{prop}[theorem]{Proposition}
\newtheorem{corol}[theorem]{Corollary}

\theoremstyle{definition}
\newtheorem{defi}[theorem]{Definition}
\newtheorem{example}[theorem]{Example}
 \newtheorem{rem}[theorem]{Remark}
 
\theoremstyle{remark}

\numberwithin{equation}{section}

\def\bd{\begin{defi}}
\def\ed{\end{defi}}
\def\bt{\begin{theorem}}
\def\et{\end{theorem}}
\def\bl{\begin{lem}}
\def\el{\end{lem}}
\def\bp{\begin{prop}}
\def\ep{\end{prop}}
\def\br{\begin{rem}}
\def\er{\end{rem}}
\def\bc{\begin{corol}}
\def\ec{\end{corol}}
\def\bex{\begin{example}}
\def\eex{\end{example}}
\def\pr{\begin{proof}}
\def\pro{\end{proof}}
\def\eqna{\begin{eqnarray*}}
\def\eqnaa{\begin{eqnarray}}
\def\deqna{\end{eqnarray*}}
\def\deqnaa{\end{eqnarray}}

\definecolor{dark}{rgb}{0.18,0.18,0.68}
\definecolor{mydark}{rgb}{0.78,0.08,0.08}
\definecolor{crew}{rgb}{0.2,0.5,0.2}
\definecolor{mmg}{rgb}{0.31,0.50,0.23}
\definecolor{dblue}{rgb}{0.01,0.01,0.44}
\definecolor{red}{rgb}{0.57,0.11,0.15}
\definecolor{cobalt}{rgb}{0.04,0.3,0.85}
\usepackage[colorlinks,citecolor=cobalt,linkcolor=cobalt,urlcolor=cobalt,pdfpagemode=UseNone,backref = page]{hyperref}

\usepackage[backref=page]{hyperref}
\renewcommand*{\backref}[1]{}
\renewcommand*{\backrefalt}[4]{%
	\ifcase #1 %
		\or        (cited on page~#2)%
	\else      (cited on pages~#2) %
	\fi}
\hypersetup{colorlinks=true,linkcolor=black, citecolor=blue,urlcolor=black}

 \input ulem.sty
 
\makeatletter
\def\subsubsection{\@startsection{subsubsection}{3}%
  \z@{.5\linespacing\@plus.7\linespacing}{.3\linespacing}%
  {\normalfont\bfseries}}
\makeatother

\begin{document}

\title[The canonical submersion of $\mc{S}$-manifolds]{The canonical submersion of $\mc{S}$-manifolds and transverse K\"ahler-Einstein structures}

\author{Ioannis Chrysikos} 
\address{Department of Mathematics and Statistics,
Faculty of Science, Masaryk University, Kotl\' a\v rsk\' a 2, 611 37 Brno, Czech Republic}
\email{chrysikos@math.muni.cz}

\maketitle

\begin{abstract}
This paper is devoted to the study of  the  holonomy properties of $(2n+s)$-dimensional $\mc{S}$-manifolds  equipped  with their characteristic 
connection. These structures generalize Sasakian geometry to higher CR-codimensions and, when viewed as geometries with parallel skew-torsion, share many holonomy features with the Sasakian case.  We show that $\mc{S}$-manifolds of arbitrary CR-codimension $s\geq 1$ provide examples of geometries with parallel
skew-torsion whose holonomy is  reducible,   indecomposable,  and of special type. We also deduce that any $\mc{S}$-manifold
admits a locally defined Riemannian submersion over a K\"ahler manifold. We   describe the corresponding curvature relations and
establish a bijective correspondence between the K\"ahler-Einstein condition on the base space
 and a generalized $\eta$-Einstein condition on the total space.
   As every $\mc{S}$-manifold comes with a characteristic foliation whose  transverse geometry is K\"ahler,  it is natural to relate the $\eta$-Einstein condition to the transverse metric, leading to a bijection between $\eta$-Einstein $\mc{S}$-manifolds and transverse K\"ahler-Einstein metrics, extending  the classical Sasakian correspondence to arbitrary CR-codimensions.
As an application  to Ricci-flat metric connections with parallel skew-torsion,  we prove that an $\mc{S}$-manifold is $\Ric^{\nabla}$-flat  if and only if it is transverse K\"ahler-Einstein with Einstein constant $\lambda=4s$.  An illustration of this characterization is presented by a Sasakian example.
 \end{abstract}
 


\pagestyle{headings}



\section{Introduction}

\subsection{Motivation}
According to a classical theorem of  \'E. Cartan,  affine metric connections on a Riemannian manifold  $(M^{n}, g)$   fall into eight classes, each distinguished by its torsion (cf. \cite{AF}).  By the fundamental theorem of Riemannian geometry,  the torsion-free case uniquely corresponds to the Levi-Civita connection $\nabla^g$.   In the presence of  torsion $T\neq 0$,  a particularly important case arises  when the tensor $g(T(-, -),-)$  is a 3-form. Such connections share the same geodesics and Killing vector fields as the Levi-Civita connection.  Moreover, the adapted connections of  many remarkable non-integrable $\Gg$-structures, such as nearly K\"ahler structures, nearly parallel $\G_2$-structures, Sasakian  and 3-Sasakian structures,   are of this type.  
Therefore, metric  connections with totally skew-symmetric torsion (briefly, skew-torsion)  often serve as 
a natural replacement of  the Levi-Civita connection,  capturing   features of the underlying non-integrable geometry
that are not visible through $\nabla^g$.   

This  remarkable scenario often  appears  in the homogeneous naturally reductive setting, where skew-torsion naturally enters the picture as the torsion of the canonical connection  (cf. \cite{Kob2}).  Indeed, it is known that the canonical connection  serves as the adapted connection for  many homogeneous (non-integrable) $\Gg$-structures, such as homogeneous nearly   K\"ahler structures and many other examples.
Connections with skew-torsion also play an important role in (type II) string theory, leading to further links with holonomy theory, since supersymmetries correspond to spinors that are parallel with respect to such connections. As a result,  metric connections with skew-torsion have attracted significant attention over the past  decades; see, for example, \cite{B89, IP01, FrIv,  Fr03, AF, CS04,  Schoe, Pu12,  AF14,   AFF15, S18,  DGP18, CGW19,    AD20,   CMS21, MS24, BFG24},   among many other works.

 A connection with totally skew-symmetric torsion $T$ is uniquely determined by the relation
\[
\nabla=\nabla^g+\frac{1}{2}T\,.
\]
In this article, a Riemannian   manifold  endowed with such a connection  will be denoted  by $(M^n, g, T)$ and called a   \textsf{geometry with skew-torsion}. 
If the torsion 3-form $T$ is  $\nabla$-parallel, then $(M^n, g, T)$ is  referred to as a  \textsf{geometry with parallel skew-torsion}. Geometries with parallel skew-torsion lie at the heart of this paper. Compact Lie groups endowed with a bi-invariant metric, non-symmetric naturally reductive spaces, and  the non-integrable structures mentioned above provide  standard examples, although many  others are known.  Geometries with parallel skew-torsion are also relevant in the Lorentzian setting (cf.  \cite{EG22}).

For metric connections with torsion  the classical de Rham decomposition theorem is {\it not} in general  valid.  In particular, reducibility of the holonomy representation does not imply a local decomposition of the manifold as a Riemannian product, in contrast to the torsion-free Riemannian case.
  Riemannian manifolds  $(M^n, g)$  admitting a metric connection $\nabla$ with non-trivial  parallel torsion having  irreducible holonomy representation were classified in   \cite{CS04}. This classification was recently revised and completed  in  \cite{MS24}.  
  
   On the other hand, the situation in the presence of reducible holonomy is considerably more involved, and only partial results are currently available.    Recently, in   \cite{CMS21}, a locally defined Riemannian submersion, called the \textsf{standard submersion}, was constructed for any geometry with parallel skew-torsion. The standard submersion provides an effective method for the realization of new examples of geometries with parallel skew-torsion and, at the same time,  enables a systematic study of their holonomy, including the significant case of reducible holonomy.   This approach has been further developed and applied to various classes of geometric structures in subsequent works; see, for instance, \cite{ADS21, MS24} and the references therein.    Note that in many cases, including  Sasakian geometry, 
   the holonomy group is generally difficult to describe.
It is therefore often useful to study the decomposition of the tangent bundle as a representation of groups larger than the holonomy group, such as the stabilizer of the torsion 3-form $T$, or any other group   lying between them. 
  This perspective leads to the notion of \textsf{canonical $\fr{g}$-splittings}, where $\fr{g}$ is an intermediate  Lie algebra, i.e., $\fr{hol}(\nabla)\subset\fr{g}\subset\fr{stab}(T)$. This notion
 was exploited in \cite{MS24}, where one of the results shows that  any canonical $\fr{g}$-splitting induces a locally defined Riemannian submersion enjoying properties analogous to those of the standard submersion.


\subsection{Outcome}
The present work lies within this line of research. 
In \cite{BC}  we   described a new class of geometries with parallel skew-torsion
based on the notion of  \textsf{metric $f$-manifolds} and, in particular, on \textsf{$\mc{S}$-manifolds}.  	Recall that metric $f$-manifolds  arise in $2n+s$ dimensions and are generalizations of almost Hermitian geometries (case $s=0$)
and almost contact metric geometries (case $s=1$).  Moreover, $\mc{S}$-manifolds are special classes of metric   $f$-manifolds that  provide a higher-dimensional analogue of Sasakian manifolds,  the latter corresponding to the case  $s=1$. 
Metric $f$-structures   have a long history, originating in the foundational works of Yano, Goldberg, and Blair  \cite{Yano63, GY70, Blair70}.  
  Such structures are important not only for their existence in both odd and even dimensions, but also because  there exist manifolds   for which an $\mc{S}$-structure seems to be the best structure one can hope to obtain. Such an example is  the Lie group $\U(2)$, although many others exist 	 (cf. \cite{TK07, DL05}).  

In this article we investigate  the holonomy properties of the characteristic connection $\nabla$ on $\mc{S}$-manifolds, introduced in \cite{BC}, and   given by $\nabla=\nabla^{g}+\frac{1}{2}\sum_{i=1}^{s}\eta_{i}\wedge\dd\eta_i$ (recall that for   $s=1$ this was first presented in \cite{FrIv}).
 We show that such  geometries with parallel skew-torsion  are reducible but indecomposable, in sharp contrast to the torsion-free Riemannian scenario (see \cite{TP13} and Remark \ref{RiemannS}).   This means that  $(M^{2n+s}, g, T)$ is {\it not}   locally a  product of two geometries with (parallel) skew-torsion, i.e., there are no manifolds $M_1,  M_2$ such that $(M^{2n+s}, g, T)=(M_1, g_1, T_1)\times (M_2, g_2, T_2)$, with $T=T_1+T_2$ and $T_i\in\Lambda^{3}TM_{i}$ for $i=1, 2$.    We further deduce that  such non-integrable $\Gg$-structures provide  new examples
 of geometries with parallel skew-torsion of \textsf{special type}.  Recall that  a geometry with parallel skew-torsion $(M^{n}, g, T)$ is said to be of  special type if the holonomy algebra 
 $\fr{hol}(\nabla)$ of the characteristic connection $\nabla$ acts trivially on the vertical distribution $\mc{V}\neq 0$ of the standard submersion, see  \cite{CMS21, MS24}.  It is known that Sasakian manifolds endowed with their unique characteristic connection $\nabla$ are   examples of such geometries and it turns out that  their higher  $(2n+s)$-dimensional analogues, i.e., $\mc{S}$-manifolds, enjoy the same property.   

Given a geometry with parallel skew-torsion $(M^n, g, T)$ one can always  consider  
a canonical $\fr{stab}(T)$-splitting $TM=\mc{H}\oplus\mc{V}$ of $TM$  into horizontal and vertical distributions, where $\fr{stab}(T)\subset\fr{so}(n)$ is the (infinitesimal) stabilizer of the torsion 3-form $T$, see \cite{MS24}. For an $\mc{S}$-manifold  viewed as a geometry with parallel skew-torsion and denoted as $(M^{2n+s}, g, T)$,   we observe that the  stabilizer of   $T$ is isomorphic to the Lie algebra $\fr{u}(n)$.  Moreover,  the splitting $TM=\mc{D}\oplus\mc{D}^{\perp}$  of the tangent bundle of $M^{2n+s}$ induced  by the $f$-structure $\phi$  can be interpreted as a canonical $\fr{u}(n)$-splitting. 
Hence, one can introduce a locally defined Riemannian submersion of  $M^{2n+s}$ over  an almost Hermitian space $N^{2n}$  with totally geodesic  fibers tangent to $\mc{D}^{\perp}$, which we prove that is    K\"ahler, see Theorem \ref{standSm}. 
We denote this map by $\pi : M^{2n+s}\to N^{2n}$ 
and refer to as the \textsf{canonical submersion} of an $\mc{S}$-manifold $M^{2n+s}$ (for   $s=1$ the canonical $\fr{u}(n)$-submersion of a Sasakian manifold was recently studied in \cite{MS24}). 
It turns out that, when $M^{2n+s}$ is a compact regular $\mc{S}$-manifold,   this   submersion is a principal torus bundle $\pi_s : M^{2n+s}\to N^{2n}$, first introduced by Blair, Ludden and Yano \cite{BLY73}.  
 
For a geometry with parallel skew-torsion $(M^{n}, g, T)$, an important situation arises
when the vertical space $\mc{V}$ of the canonical $\fr{stab}(T)$-splitting is one-dimensional and the horizontal space $\mc{H}$ is an irreducible representation of $\fr{stab}(T)$. 
In \cite{MS24}  this situation is   referred to as the \textsf{almost irreducible case}.  We may extend this to allow $\mc{V}$ to be $s$-dimensional, a situation that we  refer to as  \textsf{holonomy $s$-irreducible}, see Definition \ref{stabTHM}.  This enables us to provide a holonomy characterization of  $\mc{S}$-manifolds and  establish a higher-dimensional analogue of the corresponding   characterization of  Sasakian manifolds,   presented in  \cite[Theorem 7.1]{MS24}. This is the content of  Theorem \ref{stabTHM}.

We further investigate the curvature properties of an $\mc{S}$-manifold $M^{2n+s}$ with respect to its characteristic connection $\nabla$ and, moreover, derive  curvature identities  of the canonical submersion $\pi : M^{2n+s}\to N^{2n}$. 
We  prove that the K\"ahler-Einstein condition on the leaf space $N^{2n}$ is equivalent to a (generalized) $\eta$-Einstein condition on the total space, see Theorem \ref{KEcor}. 

 As every $\mc{S}$-manifold comes with a characteristic foliation $\scF$ whose transverse geometry is K\"ahler (cf.  \cite{DKWP}),  
 it is natural to 
 relate this correspondence  with the transverse metric. We show that there is a bijective correspondence between $\eta$-Einstein $\mathcal{S}$-structures and transverse K\"ahler-Einstein metrics, thereby extending the classical correspondence from the Sasakian setting  to the higher-dimensional context considered here; see \cite[Theorem 11.1.3]{BG08} and   \cite{FOW09} for  the  Sasakian case.  
 We finally prove a result concerning Ricci-flat metric connections with (parallel) skew-torsion, in particular we show that an $\mc{S}$-manifold is $\Ric^{\nabla}$-flat with respect to its characteristic connection $\nabla$ if and only if it is transverse K\"ahler-Einstein with Einstein constant $\lambda=4s$. Therefore we obtain a characterization of Ricci-flat metric connections with (parallel) skew-torsion that is applicable in both odd and even dimensions, since $\mc{S}$-manifolds  may occur in arbitrary CR-codimensions.   We illustrate this characterization in terms of the canonical Sasakian structure on the 3-sphere.

The structure of the paper is as follows.  In the remainder of this introductory section we collect the necessary preliminaries. 
 Section \ref{append1} is devoted to the curvature of the transverse K\"ahler metric associated  to the characteristic foliation defined on 
an $\mc{S}$-manifold.  The inclusion of this material at this point is motivated by clarity and the convenience of the reader, since many of the results established here will be needed at several points in the subsequent sections. These results are also of independent interest.    Section~\ref{CSmnfds} is devoted to the curvature properties  of the characteristic connection on  
arbitrary $\mc{S}$-manifolds, including explicit formulas  that relate the Ricci  tensor and the scalar curvature with the corresponding Riemannian versions. It also contains  the characterization of $\Ric^{\nabla}$-flat $\mc{S}$-manifolds, 
 together with examples illustrating many of the main results.  Section~\ref{SubmS}  focuses  on the holonomy properties of $\mc{S}$-manifolds endowed with their characteristic connection. Here we introduce the canonical submersion, describe the corresponding curvature identities and establish the correspondence between the generalized $\eta$-Einstein condition  on the total space and the K\"ahler-Einstein condition on the leaf space.

 %
\bigskip
\noindent\textbf{Acknowledgements.}
The author  is grateful to Anton Galaev (UHK) and Paul Schwahn (Unicamp)  for many helpful discussions and comments that improved the manuscript. He acknowledges the support of the Czech Science Foundation (project GA24-10887S) and the Horizon 2020 MSCA project CaLIGOLA (ID 101086123).


\subsection{Preliminaries on metric $f$-manifolds}\label{Preliminaries} 
   Globally framed $f$-manifolds are $(2n+s)$-dimensional  manifolds
endowed with an \textsf{$f$-structure} $\phi : TM\to TM$ of constant rank $2n$, such that 
the distribution $\Ker(\phi)\subset TM$ is parallelizable. Thus, 
 there exist global vector fields $\xi_1, \ldots, \xi_s\in TM$ on $M$,  called the \textsf{characteristic vector fields}, and  smooth 1-forms $\eta_1, \ldots, \eta_s\in T^*M$  satisfying 
\[
\Ker(\phi)=\langle \xi_1, \ldots, \xi_s\rangle\,,\quad   \eta_{i}(\xi_j)=\delta_{ij}\,, \quad \eta_{i}\circ\phi=0\,,\quad  \phi^2=-\Id+\sum_{i=1}^{s}\eta_i\otimes\xi_i\,,
\]
for all $i, j\in\{1, \ldots, s\}$, see  \cite{Blair70}.   In this case  the tangent bundle  $TM$ of $M$ splits into two complementary subbundles $\Imm(\phi)$ and $\Ker(\phi)$, which we denote by $\mc{D}, \mc{D}^{\perp}$, respectively, 
  \[
  TM=\mc{D}\oplus\mc{D}^{\perp}\,,\quad \mc{D}:=\Imm(\phi)=\bigcap_{i=1}^{s}\Ker(\eta_i)\,,\quad \mc{D}^{\perp}:=\Ker(\phi)=\langle \xi_1, \ldots, \xi_s\rangle\,.
  \]
  Tangent vectors belonging to the distribution $\mc{D}$  are called \textsf{horizontal} and tangent vectors belonging to $\mc{D}^{\perp}$ are called \textsf{vertical}.  
  
For any globally framed $f$-structure $(\phi, \xi_i, \eta_j)$   there exists an adapted Riemannian metric,  i.e., a Riemannian metric $g$ such that
\begin{equation}\label{gXgY}
g(\phi X, \phi Y)=g(X, Y)-\sum_{i=1}^{s} \eta_{i}(X)\eta_{i}(Y)\,,\quad   \ X, Y\in TM\,.
 \end{equation}
 Then   the family $(\phi, \xi_i, \eta_j, g)$ is   referred to as a \textsf{metric $f$-structure} and  $(M^{2n+s}, \phi, \xi_i, \eta_j, g)$ is said to be a \textsf{metric $f$-manifold} (cf. \cite{Blair70, CFF90} and for more details we refer to \cite{BC}). 
Note that the horizontal distribution $\mc{D}$ is a complex vector bundle with almost complex structure defined by  the restriction $\phi|_{\mc{D}}$, since   $\phi^2|_{\mc{D}}=-\Id_{\mc{D}}$. 
Thus,   $(M^{2n+s}, \phi, \xi_i, \eta_j, g)$   naturally exhibits the structure of an almost CR-manifold of CR-dimension $n$ and  CR-codimension $s=\dim\Ker(\phi)$. Obviously, the case $s=0$ corresponds to almost Hermitian manifolds, while for $s=1$ we obtain the notion of almost contact metric manifolds.  Therefore, metric $f$-manifolds 
 should  be understood as higher-dimensional generalizations of almost Hermitian manifolds and almost contact metric manifolds.
 
Let  $(M^{2n+s}, \phi, \xi_i, \eta_j, g)$ be a metric $f$-manifold. From (\ref{gXgY}) it follows that
\begin{equation}\label{gphi}
g(\phi X, Y)+g(X, \phi Y)=0\,, 
\end{equation}
  for all $X, Y\in TM$. Hence the rule $F(X, Y):=g(X, \phi Y)$ defines  a 2-form  $F$ on $M^{2n+s}$, 
  called  the \textsf{fundamental 2-form}. We also recall the Nijenhuis tensor  $N^{(1)} : TM\times TM\to TM$, 
\[
 N^{(1)}:=N_{\phi}+\sum_{i=1}^{s}\dd\eta_i\otimes\xi_i\,,\quad 
 N_{\phi}(X, Y):=[\phi X, \phi Y]+\phi^{2}[X, Y]-\phi[X, \phi Y]-\phi[\phi X, Y]\,,
\]
  for all $X, Y\in TM$.  
\bd
A metric $f$-manifold $(M^{2n+s}, \phi, \xi_i, \eta_j, g)$ is said to be
\begin{itemize}
\item \textsf{normal} if $N^{(1)}=0$;
\item \textsf{contact metric $f$-manifold} or \textsf{almost $\mc{S}$-manifold} if $s\geq 1$ and $2F=\dd\eta_j$ for all $j\in\{1, \ldots, s\}$;
\item \textsf{$\mc{K}$-manifold}, if $\dd F=0$ and $N^{(1)}=0$;
\item \textsf{$\mc{S}$-manifold} if $s\geq 1$, $2F=\dd\eta_j$ for all $j\in\{1, \ldots, s\}$ and $N^{(1)}=0$.
\end{itemize}
\ed

For any normal metric $f$-manifold of CR-codimension $s\geq 1$ we have (cf. \cite{GY70})
 \[
\mc{L}_{\xi_i}\xi_j=[\xi_i, \xi_j]=0\,,\quad \mc{L}_{\xi_i}\eta_{j}=0=\mc{L}_{\xi_i}\phi\,,\quad \dd\eta_{i}(\phi X, Y)+\dd\eta_{i}(X, \phi Y)=0\,,
\]
 for all vector fields $X, Y\in TM$ and $i, j\in\{1, \ldots, s\}$.
A $\mc{K}$-manifold for $s=0$ is a K\"ahler manifold  and  for $s=1$ a quasi-Sasakian manifold $(\dd F=0, N^{(1)}=0)$. 
Observe also that $\mc{S}$-manifolds are contact metric $f$-manifolds which are normal, and special examples
of $\mc{K}$-manifolds. Obviously,  for $s=1$ almost $\mc{S}$-manifolds are contact metric manifolds and $\mc{S}$-manifolds are Sasakian manifolds.   Note however, that   unlike Sasakian manifolds, no $\mc{S}$-structure with $s\geq 2$ can be realized on a simply connected compact manifold (cf. \cite{DL05}).    By \cite{Blair70} it is also known that  on any $\mc{K}$-manifold $(M^{2n+s},  \phi, \xi_i, \eta_j, g)$ with $s\geq 1$ the characteristic vector fields are Killing. Moreover, we have $[\xi_i, \xi_j]=0$ for all $i, j\in\{1, \ldots, s\}$
as a  consequence of the normality condition.  The characteristic vector fields commute also on   almost $\mc{S}$-manifolds (cf. \cite{CFF90}).

\subsection{The characteristic foliation on $\mc{S}$-manifolds}
Let  $(M^{2n+s}, \phi, \xi_i, \eta_j, g)$ be a   $\mc{K}$-manifold and let   $TM=\Imm(\phi)\oplus\Ker(\phi)=\mc{D}\oplus\mc{D}^{\perp}$ 
be decomposition of its tangent bundle into horizontal and vertical parts, as described in Section \ref{Preliminaries}. 
Clearly, the vertical distribution $\mc{D}^{\perp}$ is  integrable and  when $M$  is an $\mc{S}$-manifold, then $\mc{D}^{\perp}$ is flat (see  \cite[Theorem 1.7]{Blair70}). 
 Moreover,  $\mc{D}^{\perp}$ is generated by  the linearly-independent Killing vectors fields $\xi_i$ for $i\in\{1, \ldots, s\}$, hence we   obtain an $s$-dimensional Riemannian foliation $\scF$ (also denoted as $\scF_{\xi_1, \ldots, \xi_s})$  defined by the infinitesimal action of $\R^s$ on $M$ by isometries, with  tangent bundle $T\scF=\mc{D}^{\perp}$. Note that since the metric $g$ is bundle-like,   the horizontal distribution $\mc{D}$ is  totally geodesic (cf. \cite[Remark 1.6]{BGM94}).

\bd\label{characteristicF}
The Riemannian foliation $\scF$ is called the  \textsf{characteristic foliation} of  $M^{2n+s}$.
\ed
For $s=1$ we obtain the characteristic foliation of a Sasakian manifold, see \cite{BG08, W17} for more details.  
  Let us now consider the  corresponding leaf space
\[
M^{2n+s}/\scF\,.
\]
When the leaves are compact,   then the leaf space is an {\it orbifold}, see for example  \cite[Theorem 1.8]{BGM94}. 
It is known that orbifolds arise naturally as the leaf spaces of certain well-behaved Riemannian foliations (e.g., in the Sasakian setting, see \cite{BGM94, BG08}). 
 
\subsection{Regular $\mc{S}$-manifolds as principal torus bundles}\label{toruS}
According to Palais \cite{P57}, to impose a smooth structure on the leaf space  $M^{2n+s}/\scF$, it is necessary to assume the regularity condition for $\scF$ and the compactness of the leaves. In this case, the map 
\[
\pi_{s} : M^{2n+s}\to M^{2n+s}/\scF
\]
 becomes a  smooth fibration having as fibers  the leaves of $\scF$ (we always assume that  $M^{2n+s}$ is connected).  When the characteristic foliation $\scF$ is regular and the characteristic vector fields $\xi_i$ are also regular,  then one  speaks for a \textsf{regular $\mc{K}$-structure} $(\phi, \xi_i, \eta_j, g)$ and \textsf{regular $\mc{K}$-manifold} $(M^{2n+s}, \phi, \xi_i, \eta_j, g)$, respectively. The definitions of `regularity' that we avoid to recall here can be found in \cite[p.~177]{BLY73},   see also the discussion in \cite[pp.~132-135]{FIP04} and   \cite[Section~3]{BP12}.
 
 For compact regular $\mc{K}$-manifolds the maximal integral curves of each characteristic vector field are diffeomorphic to the circle $\Ss^1$ and hence
 the leaves of $\scF$ (or, equivalently, the fibers of $\pi_{s}$)
  are diffeomorphic  to the $s$-dimensional torus $\Tg^s$
generated by the flow of     $\xi_i$'s,
\[
\Tg^s=\overline{\langle\exp(t_i\xi_i), \ldots, \exp(t_s\xi_s)\rangle}\,.
\]
This is contained in the isometry group of $(M^{2n+s}, \phi, \xi_i, \eta_j, g)$ since it is generated by Killing vector fields (cf. \cite{Bes}). 
We thus obtain a principal torus bundle
 \[
 \Tg^s\longrightarrow M^{2n+s}\overset{\pi_s}{\longrightarrow} M^{2n+s}/\scF\,.
 \]
 \bt \textnormal{(\cite{BLY73, S81})}\label{BYL}
Let $(M^{2n+s}, \phi, \xi_i, \eta_j, g)$ be a  compact  regular $\mc{K}$-manifold. Then, for any $s\geq 1$,  $M^{2n+s}$ is the total space of a principal $\Tg^s$-bundle
\[
\pi_s : M^{2n+s}\to N^{2n}=M^{2n+s}/\scF
\]
 over a  $2n$-dimensional compact K\"ahler manifold $N=M/\scF$.  In particular,  $\pi_s$ is a Riemannian submersion with totally geodesic fibers. If in addition $M^{2n+s}$ is an $\mc{S}$-manifold, then $N^{2n}$ is a Hodge manifold for any $s\geq 1$, i.e., $[\Om]\in H^{2}(N;\Z)$.
\et
This result generalizes the corresponding statement in the regular Sasakian case $(s=1)$, 
see \cite[Theorem~7.5.1, (iv)]{BG08}.   
Finally, we recall a construction concerning the converse direction.
\bt\label{BS70}
{\rm(1)} \textnormal{(\cite{Blair70})} Let $N^{2n}$ be a  K\"ahler manifold and let  $M^{2n+s}$ be the total space of a principal torus bundle $\pi : M^{2n+s}\to N^{2n}$ over  $N$. Let $\gamma=(\eta_1, \ldots, \eta_s)$ be the  Lie algebra valued connection 1-form on $M^{2n+s}$ such that $\dd\eta_i=\pi^*\Om$, for any $i\in\{1, \ldots, s\}$, where $\Om$ is the K\"ahler form on $N^{2n}$. Then $M^{2n+s}$ is an $\mc{S}$-manifold.\\
{\rm(2)}  \textnormal{(\cite{S81})}  If  in addition $N^{2n}$ is assumed to be a Hodge manifold, i.e., when the K\"ahler form $\Om$ is an integral (symplectic) 2-form on $N$,  then for any $s\geq 1$ there exists a principal torus bundle  $\pi : M^{2n+s}\to N^{2n}$ over $N^{2n}$ whose total space $M^{2n+s}$  has a regular $\mc{S}$-structure.
\et

Next, for simplicity,  we will adopt the following definition.
\bd
 The Riemannian submersion $\pi_s : M^{2n+s}\to N^{2n}$ of  a compact regular $\mc{S}$-manifold $M^{2n+s}$  is called the \textsf{$\BLY$-submersion}. 
 \ed
 We will denote by $\scH^{2n}, \scV^{s}$ the horizontal and vertical distributions of $\pi_s$, such that $TM=\scH^{2n}\oplus\scV^s$. 
Obviously,  by the definition of $\pi_s$ we have  
\begin{equation}\label{vh}
\scH_{x}^{2n}=\mc{D}_{x}=\Imm(\phi_x)=\bigcap_{i=1}^{s}\Ker((\eta_i)_{x})\,,\quad \scV_x^{s}=\mc{D}^{\perp}_{x}=\Ker(\phi_x)=\langle (\xi_1)_{x}, \ldots, (\xi_s)_{x}\rangle
\end{equation}
 at any $x\in M$, hence we identify $\scH^{2n}=\mc{D}$ and $\scV^{s}=\mc{D}^{\perp}$, respectively.

\br Given a $\BLY$-submersion $\pi_s : M^{2n+s}\to N^{2n}$ with $s\geq 2$,   one can construct 
a further submersion $\tau_{s-1} : M^{2n+s}\to\Xg^{2n+1}$ over the $(2n+1)$-dimensional compact regular Sasakian manifold associated to the K\"ahler manifold $N$ by the $\BLY$-submersion for $s=1$. This makes the following diagram
commutative.
\[
  \xymatrix{
  M^{2n+s} \ar[rr]^{\tau_{s-1}} \ar[dr]_{\pi_s}  & & \Xg^{2n+1}  \ar[dl]^{\pi_1}  \\
& N^{2n}& }        
\]
Note that $\pi_1$ is a principal $\Ss^1$-fibration, while $\tau_{s-1}$ is a principal torus bundle over $\Xg^{2n+1}$, 
with fiber diffeomorphic to the $(s-1)$-dimensional torus $\Tg^{s-1}$.  Moreover, all the maps in this diagram are Riemannian submersions with totally geodesic fibers.
\er


 \section{The transverse K\"ahler geometry of the characteristic foliation}\label{append1}

 \subsection{Curvature identities of the  transverse K\"ahler metric}

   This section deals with the transverse geometry of the   characteristic foliation $\mc{F}$ on a $\mc{S}$-manifold  $(M^{2n+s}, \phi, \xi_i, \eta_j, g)$ defined by the characteristic   vector fields $\xi_1, \ldots, \xi_s$.  
   We establish curvature identities that extend to arbitrary CR-codimesnion $s\geq 1$ the classical identities characterising the transverse geometry of  Sasakian manifolds ($s=1$), see 
\cite{BGM96, BG08, FOW09} for further details.

Many features of the transverse K\"ahler geometry considered here were studied in \cite{DKWP}, primarily in the context of $\mc{K}$-structures; see also \cite[Section 5]{TK07} and \cite{GL20,R23} for further results on the transverse geometry of the foliation determined by the kernel $\Ker(\phi)$ of various metric $f$-structures $(\phi,\xi_i,\eta_j,g)$. 
However, for $s>1$,  the curvature identities derived below do not appear to have been recorded in the literature.
As a number of the results presented below, in particular the formulas for the \textsf{transverse Ricci tensor}, will be used in the forthcoming sections, we include this exposition here for the reader's convenience.
For an introduction to foliations and basic tensor fields, we refer  to \cite{T97}; see also \cite{BG08,AC22,R23}.

Let us fix, once and for all,   an $\mc{S}$-manifold  $(M^{2n+s}, \phi, \xi_i, \eta_j, g)$.  Recall that in this case we have $[\xi_i, \xi_j]=0$ for all $i, j\in\{1, \ldots, s\}$, each characteristic vector field $\xi_i$ is  Killing and  the following relations hold (see also \cite{CFF90})
\[
2F=\dd\eta_1=\cdots=\dd\eta_s\,,\quad N_{\phi}=-\sum_{i=1}^{s}\dd\eta_i\otimes\xi_i\,,\quad \nabla^{g}_{X}\xi_i=-\phi(X)\,,\quad \forall \ X\in TM\,.
\]
Moreover, $[X, \xi_i]\in\mc{D}$ for all $X\in\mc{D}$. 
\bl\label{XYvert}
Let $(M^{2n+s}, \phi,\xi_i, \eta_j, g)$ be an $\mc{S}$-manifold. Then, for any two horizontal vector fields $X, Y\in\mc{D}$ we have
\[
[X, Y]_{\mc{D}^{\perp}}=-2\sum_{i=1}^{s}F(X, Y)\xi_i\,.
\]
\el
\pr
By definition,  $\phi(W)=0$ for any $W\in\mc{D}^{\perp}$, where $\mc{D}^{\perp}=\langle\xi_1, \ldots, \xi_s\rangle$.  Thus $\phi^2(W)=0$,  that is,  $W=\sum_{i}\eta_{i}(W)\xi_{i}$. 
Moreover, $[X, Y]_{\mc{D}^{\perp}}\in\mc{D}^{\perp}$  and $\mc{D}=\bigcap_{i}\Ker(\eta_i)$, we thus  get
\[
[X, Y]_{\mc{D}^{\perp}}=\sum_{i}\eta_{i}([X, Y]_{\mc{D}^{\perp}})\xi_i=\sum_{i}\eta_{i}([X, Y])\xi_i=-\sum_{i}\dd\eta_{i}(X, Y)\xi_i=-2\sum_{i}F(X, Y)\xi_i\,,
\]
since  for any $X, Y\in\mc{D}$ and $i\in\{1, \ldots, s\}$  it holds that $\dd\eta_{i}(X, Y)=-\eta_{i}([X, Y])$. 
\pro

Let $\mc{F}$ be the characteristic foliation on $M^{2n+s}$ (see Definition \ref{characteristicF}). 
A $p$-form $\al$  on $M$ is   {\textsf{basic}} if and only if it satisfies the relations
\[
\xi_i\lrcorner\al=0\,,\quad \mc{L}_{\xi_i}\al=0\,,
\]
for all $i\in\{1, \ldots, s\}$, where $\mc{L}$ denotes the Lie derivative. 
Now, by definition,  an $\mc{S}$-manifold is a $\mc{K}$-manifold; thus its characteristic foliation  $\scF$ is   a transverse Hermitian foliation which is  K\"ahler, see  \cite[Theorem 1]{DKWP}. 
In particular,   we have an exact sequence of vector bundles
\[
0\longrightarrow\mc{D}^{\perp}\longrightarrow  TM \longrightarrow TM/\mc{D}^{\perp}\longrightarrow 0
\]
where the normal bundle $TM/\mc{D}^{\perp}$ has a transverse complex structure $\check{J}$ and a compatible transverse Riemannian metric $\check{g}$ such that  the 2-form $\check{\omega}$ induced by $\check{J}$ and $\check{g}$ is basic and symplectic. This is the horizontal  2-form 
\[
\check{\om}:=F|_{\mc{D}\times\mc{D}}=\frac{1}{2}\dd\eta_j|_{\mc{D}\times\mc{D}}
\]
for some (and thus any)  $j\in\{1, \ldots, s\}$, which is a symplectic form on the horizontal distribution $\mc{D}=\bigcap_{j=1}^{s}\Ker(\eta_j)$. It is basic since $\xi_i\lrcorner F=0$ and $\xi_i\lrcorner \dd F=0$ for all $i\in\{1, \ldots, s\}$ (see also \cite[p.~156]{Blair70}). Here, we have identified the complex vector bundles $(TM/\mc{D}^{\perp}, \check{J})$ and  $(\mc{D}, \phi|_{\mc{D}})$, and in particular   $\check{J}$ with $\phi|_{\mc{D}}$.  The transverse mertic $\check{g}$ is then defined by 
\begin{equation}\label{TKahlerg}
\check{g}(X, Y):=\check{\omega}(\check J X, Y)=\frac{1}{2}\dd\eta_{j}(\phi X, Y)=F(\phi X, Y)=g(\phi X, \phi Y)=g(X, Y)\,,
\end{equation}
for all $X, Y\in\mc{D}$, where $j\in\{1, \ldots, s\}$.  
Thus we have $g(X, Y)=\check{g}(X,Y)+\sum_{i=1}^{s}\eta_{i}(X)\eta_{i}(Y)$, for all $X, Y\in TM$, that is,
\[
g=\check{g}+\sum_{i=1}^{s}\eta_{i}\otimes\eta_{i}\,.
\]
For $s=1$ this  yields the   transverse K\"ahler metric associated to the characteristic foliation (Reeb foliation)  $\mc{F}=\mc{F}_{\xi}$ defined on a Sasakian manifold, see \cite{BG08}.    
 Next we will describe curvature  identities characterizing the transverse  geometry of the characteristic foliation defined on an $\mc{S}$-manifold. For the transverse metric $\check{g}$ we define a connection on the horizontal distribution $\mc{D}$ by the rule
(see for example \cite[p.~21]{T97} or \cite[p.~188]{AC22})
\[
\check\nabla_{X}Y:=
\begin{cases}
(\nabla^{g}_{X}Y)_{\mc{D}}\,, & \text{if}  \ X, Y\in\mc{D}\,,\\
[X, Y]_{\mc{D}}\,, & \text{if} \ X\in\mc{D}^{\perp}, Y\in\mc{D}\,.\\
\end{cases}
\]
In particular, we have $\check{\nabla}_{\xi_i}Y:=[\xi_i, Y]_{\mc{D}}=[\xi_i, Y]$ for any $i\in\{1, \ldots, s\}$ and $Y\in\mc{D}$ (recall
that   $\mc{D}^{\perp}=\langle\xi_1, \ldots, \xi_s\rangle$
 and $[X, \xi_i]\in\mc{D}$ for any  $i\in\{1, \ldots, s\}$). 
One can easily check that $\check{\nabla}$ is the unique torsion-free  connection preserving $\check{g}$, thus we will refer to $\check{\nabla}$
as the \textsf{transverse Levi-Civita connection}. This is  also  known  as the  \textsf{(adapted) Bott connection} on $\mc{D}$ (cf. \cite[p.~135]{DTPW14}).  
 \bl\label{transvLC} Let $(M^{2n+s}, \phi, \xi_i, \eta_j, g)$ be an $\mc{S}$-manifold.  Then,  for any $X, Y\in\mc{D}$ and $i\in\{1, \ldots, s\}$
 the following hold:
  \[
 \nabla^{g}_{\xi_i}X=-\phi(X)+[\xi_i, X]=\check\nabla_{\xi_i}X-\phi(X)\,,\quad
   \nabla^{g}_{X}Y=\check\nabla_{X}Y-\sum_{i=1}^{s}F(X, Y)\xi_i\,.
 \]
 \el
 \pr
 The first relation  relies on the torsion-free property of $\nabla^g$, the relation  $\nabla^{g}_{X}\xi_i=-\phi(X)$ for all $i$ and $X\in\mc{D}$ and the definition of $\check\nabla$.   For the second  one, we   express the vector field $\nabla^{g}_{X}Y$ as 
 \[
 \nabla^{g}_{X}Y=(\nabla^{g}_{X}Y)_{\mc{D}}+(\nabla^{g}_{X}Y)_{\mc{D}^{\perp}}\,.
 \]
For the vertical projection,  since $\mc{D}^{\perp}=\langle\xi_1, \ldots, \xi_s\rangle$ and $g(\xi_i, \xi_j)=\delta_{ij}$, we  deduce that
 \begin{eqnarray*}
 (\nabla^{g}_XY)_{\mc{D}^{\perp}}&=&\sum_{i=1}^{s}g(\nabla^{g}_{X}Y, \xi_i)\xi_i=-\sum_{i=1}^{s}g(\nabla^{g}_{X}\xi_i, Y)\xi_i=\sum_{i=1}^{s}g(\phi(X), Y)\xi_i\\
 &\overset{(\ref{gphi})}{=}&-\sum_{i=1}^{s}g(X, \phi(Y))\xi_i=-\sum_{i=1}^{s}F(X, Y)\xi_i\,,
 \end{eqnarray*}
 where we employed  the metric property of $\nabla^g$ to obtain the second equality. 
 Therefore, the given expression follows  in a combination with the definition of $\check\nabla$.
 \pro
 
  \bp\label{transvJ}
 Let $(M^{2n+s}, \phi, \xi_i, \eta_j, g)$ be an $\mc{S}$-manifold. Then 
 $\check{\nabla}\check{J}=0$ and hence also $\check\nabla F=0$.
 \ep
 \pr
  Let $X, Y\in\mc{D}$. We use the decomposition  $\nabla^{g}_{X}Y=(\nabla^{g}_{X}Y)_{\mc{D}}+(\nabla^{g}_{X}Y)_{\mc{D}^{\perp}}$ and the relation $\mc{D}^{\perp}=\Ker(\phi)$ to get $\phi(\nabla^{g}_{X}Y)=\phi\big((\nabla^{g}_{X}Y)_{\mc{D}}\big)$. Thus
\begin{eqnarray*}
(\check{\nabla}_{X}\check{J})Y&=&\check\nabla_{X}\phi(Y)-\phi(\check\nabla_{X}Y)\\
&=&\big(\nabla^{g}_{X}\phi(Y)\big)_{\mc{D}}-\phi\big((\nabla^{g}_{X}Y)_{\mc{D}}\big)\\
&=&\big(\nabla^{g}_{X}\phi(Y)-\phi(\nabla^g_{X}Y)\big)_{\mc{D}}=\big((\nabla^g_{X}\phi)Y\big)_{\mc{D}}\,.
\end{eqnarray*}
However, it is known that
\[
g((\nabla^{g}_{X}\phi)Y, Z)=\sum_{i}\big\{\eta_{i}(Z)F(\phi Y, X)+\eta_{i}(Y)F(X, \phi Z)\big\}\,, \quad X, Y, Z\in TM\,,
\]
 see \cite[Lemma.~2.2]{CFF90}.  
 A combination of these two relations gives $g\big(\big((\nabla^g_{X}\phi)Y\big)_{\mc{D}}, Z\big)=0$ for all $Z\in\mc{D}$ and thus for all $X, Y\in\mc{D}$ we get  $(\check{\nabla}_{X}\check{J})Y=0$.  Moreover, we see that
 \begin{eqnarray*}
 (\check\nabla_{\xi_i}\check J)Y&=&\check\nabla_{\xi_i}\phi(Y)-\phi(\check\nabla_{\xi_i}Y)
 =[\xi_i, \phi(Y)]_{\mc{D}}-\phi([\xi_i, Y]_{\mc{D}})\\
 &=&[\xi_i, \phi(Y)]-\phi([\xi_i, Y])=(\mc{L}_{\xi_i}\phi)Y=0\,.
 \end{eqnarray*}
 Here we are based on the fact that $M^{2n+s}$ is normal, so we have  $\mc{L}_{\xi_i}\phi=0$ for all $i\in\{1, \ldots, s\}$. 
 This proves that $\check\nabla\check J=0$. Combining this relation with the $\check\nabla$-parallelism of the transverse metric, i.e., $\check\nabla\check{g}=0$, we also get $\check\nabla F=0$.  
  \pro
The  curvature tensor    induced by the transverse Levi-Civita connection $\check{\nabla}$ is defined by
\[
\check R(X, Y)Z:=\check\nabla_{X}\check\nabla_{Y}Z-\check\nabla_{Y}\check\nabla_{X}Z-\check\nabla_{[X, Y]}Z\,.
\]
It is therefore natural to refer to this curvature tensor as the  \textsf{transverse curvature tensor}. At this point it will be useful to recall known   identities for the Riemannian curvature tensor  and the Riemannian Ricci tensor on  an $\mc{S}$-manifold.  
 \bp\textnormal{(\cite[Proposition~3.4]{CFF90})}\label{known1}
Let  $(M^{2n+s}, \phi,\xi_i, \eta_j, g)$   be an $\mc{S}$-manifold. Then the Riemannian curvature tensor $R^{g}$ satisfies the following identities:\\
\[
R^{g}(X, Y)\xi_i=\sum_{j=1}^{s}\left\{\eta_{j}(X)\phi^2Y-\eta_{j}(Y)\phi^2X\right\}\,,\qquad R^{g}(X, \xi_i)Y=-(\nabla^g_{X}\phi)Y
\]
for all $X, Y\in TM$ and $i\in\{1, \ldots, s\}$.  Moreover, the Riemannian Ricci tensor 
satisfies 
 \begin{equation}\label{not_diag}
 \Ric^{g}(X, \xi_\ell)=2n\sum_{j=1}^{s}\eta_{j}(X)\,,\quad \forall \ X\in TM,\ \ell\in\{1, \ldots, s\}\,.
 \end{equation}
\ep
By  (\ref{not_diag}) we deduce that for $s>1$ the Ricci tensor $\Ric^{g}$ cannot be diagonal, and in particular an $\mc{S}$-manifold of CR-codimension $s>1$ is never Einstein, see also \cite{KT72}. 
For example,  already for $s=2$ it appears the non-zero component $\Ric^{g}(\xi_1, \xi_2)=2n$, where we assume that $\mc{D}^{\perp}=\langle\xi_1, \xi_2\rangle$.

\bp
Let $(M^{2n+s}, \phi, \xi_i, \eta_j, g)$ be an $\mc{S}$-manifold.  Then,  the transverse   curvature tensor satisfies  $\check{R}(X, \xi_i)Y=0$ for all $i\in\{1, \ldots, s\}$ and $X, Y\in\mc{D}$.
Moreover, 
\begin{equation}\label{transvR}
R^{g}(X, Y)Z=\check{R}(X, Y)Z+s\big\{F(Y, Z)\phi(X)-F(X, Z)\phi(Y)\big\}-2s\,F(X, Y)\phi(Z)\,,
\end{equation}
for any $X, Y, Z\in\mc{D}$. Therefore, 
\begin{equation}\label{transvR4}
R^{g}(X, Y, Z, W)=\check{R}(X, Y, Z, W)+s\big\{F(X, Z)F(Y, W)-F(Y, Z)F(X, W)\big\}+2s\,F(X, Y)F(Z, W)\,,
\end{equation}
for all $X, Y, Z, W\in\mc{D}$, where  we set $\check{R}(X, Y, Z, W):=\check{g}(\check R(X, Y)Z, W)$. 
 \ep

\pr
 A direct computation based on the definition of $\check\nabla$ shows that
\[
\check{R}(X, \xi_i)Y=\Big(R^{g}(X, \xi_i)Y+(\nabla^{g}_{X}\phi)Y\Big)_{\mc{D}}
\]
and the first claim follows by Proposition \ref{known1}. \\
Consider now   arbitrary horizontal vector fields $X, Y, Z\in\mc{D}$ and set $\bar{\xi}:=\sum_{i=1}^{s}\xi_i$.  By Lemma \ref{transvLC} and  the Leibniz rule we see that
\begin{eqnarray*}
\nabla^{g}_{X}\nabla^{g}_{Y}Z&=&\nabla^{g}_{X}\big(\check\nabla_{Y}Z-F(Y, Z)\bar{\xi}\big)=\nabla^{g}_{X}\check{\nabla}_{Y}Z-X\big(F(Y, Z)\big)\bar\xi-F(Y, Z)\sum_{i=1}^{s}\nabla^{g}_{X}\xi_i\\
&=&\check{\nabla}_{X}\check\nabla_{Y}Z-F(X, \check{\nabla}_{Y}Z)\bar{\xi}-X\big(F(Y, Z)\big)\bar\xi+s\,F(Y, Z)\phi(X)\,,
\end{eqnarray*}
and similarly
\begin{eqnarray*}
\nabla^{g}_{Y}\nabla^{g}_{X}Z&=&\check{\nabla}_{Y}\check\nabla_{X}Z-F(Y, \check{\nabla}_{X}Z)\bar{\xi}-Y\big(F(X, Z)\big)\bar\xi+s\,F(X, Z)\phi(Y)\,.
\end{eqnarray*}
 Recall now that the Lie bracket $[X, Y]$ of two horizontal vector fields is not necessarily horizontal;  by  Lemma \ref{XYvert} we have
\[
[X, Y]=[X, Y]_{\mc{D}}+[X, Y]_{\mc{D}^{\perp}}=[X, Y]_{\mc{D}}-2F(X, Y)\sum_{i=1}^{s}\xi_i=[X, Y]_{\mc{D}}-2F(X, Y)\bar\xi\,,\quad\forall \ X, Y\in\mc{D}\,.
\]
Thus, by Lemma \ref{transvLC} we obtain that
\begin{eqnarray*}
\nabla^{g}_{[X, Y]}Z&=&\nabla^{g}_{[X, Y]_{\mc{D}}}Z+\nabla^{g}_{[X, Y]_{\mc{D}^{\perp}}}Z\\
&=&\check\nabla_{[X, Y]_{\mc{D}}}Z-F([X, Y]_{\mc{D}}, Z)\bar{\xi}+\nabla^{g}_{-2F(X, Y)\bar\xi}Z\\
&=&\check\nabla_{[X, Y]_{\mc{D}}}Z-F([X, Y]_{\mc{D}}, Z)\bar{\xi}-2F(X, Y)\sum_{i=1}^{s}\nabla^{g}_{\xi_i}Z\\
&=&\check\nabla_{[X, Y]_{\mc{D}}}Z-F([X, Y]_{\mc{D}}, Z)\bar{\xi}-2F(X, Y)\sum_{i=1}^{s}\big\{[\xi_i, Z]-\phi(Z)\big\}\\
&=&\Big\{\check\nabla_{[X, Y]_{\mc{D}}}Z-2F(X, Y)\sum_{i=1}^{s}[\xi_i, Z]\Big\}-F([X, Y]_{\mc{D}}, Z)\bar{\xi}+2sF(X, Y)\phi(Z)\\
&=&\check\nabla_{[X, Y]}Z-F([X, Y]_{\mc{D}}, Z)\bar{\xi}+2sF(X, Y)\phi(Z)\,,
\end{eqnarray*}
since by the definition of $\check\nabla$ it is easy to see that 
\[
\check\nabla_{[X, Y]}Z=\check\nabla_{[X, Y]_{\mc{D}}}Z-2F(X, Y)\sum_{i=1}^{s}[\xi_i, Z]\,.
\]
 Altogether, 
this gives
\[
R^{g}(X, Y)Z=\check{R}(X, Y)Z+s\big\{F(Y, Z)\phi(X)-F(X, Z)\phi(Y)\big\}-2sF(X, Y)\phi(Z)+{\sf H}(X, Y, Z)\bar{\xi}
\]
where ${\sf H}(X, Y, Z)$ is  defined by
\[
{\sf H}(X, Y, Z):=-F(X, \check\nabla_{Y}Z)+F(Y, \check\nabla_{X}Z)-X\big(F(Y, Z)\big)+Y\big(F(X, Z)\big)+F([X, Y]_{\mc{D}}, Z)\,.
\]
We will show that ${\sf H}(X, Y, Z)=0$ for all $X, Y, Z\in\mc{D}$. We rely on the $\check\nabla$-parallelism of $F$, proved in Proposition \ref{transvJ}.  In particular, the condition $\check\nabla F=0$ gives that
\[
X\big(F(Y, Z)\big)=F(\check\nabla_{X}Y, Z)+F(Y, \check\nabla_{X}Z)\,,\quad
Y\big(F(X, Z)\big)=F(\check\nabla_{Y}X, Z)+F(X, \check\nabla_{Y}Z)\,.
\]
Thus,  
\begin{eqnarray*}
{\sf H}(X, Y, Z)&=&-F(\check\nabla_{X}Y, Z)+F(\check\nabla_{Y}X, Z)+F([X, Y]_{\mc{D}}, Z)\\
&=&-F([X, Y]_{\mc{D}}, Z)+F([X, Y]_{\mc{D}}, Z)=0\,,
\end{eqnarray*}
where we used that $\check\nabla$ is torsion-free. The expression in 
(\ref{transvR4}) follows immediately from  (\ref{transvR}), (\ref{gphi}), and the definition of $F$.
\pro

\br
For $s=1$ the expression (\ref{transvR4}) reduces to the formula given in part (iii) of \cite[Lemma 7.3.10]{BG08} for $(2n+1)$-dimensional $K$-manifolds (and hence also for Sasakian manifolds).   
   We mention that our approach  differs  from that in \cite{BG08} for $s=1$, as in this part we do not rely on the O'Neill's tensors.  
 \er
 \subsection{The transverse Ricci tensor}
Given an    $\mc{S}$-manifold  $(M^{2n+s}, \phi, \xi_i, \eta_j, g)$, the  \textsf{transverse Ricci tensor}  is the Ricci tensor associated to the transverse Levi-Civita connection, i.e.,
\[
\Ric^{\check{g}}(X, Y):=\sum_{a=1}^{2n}\check{g}(\check{R}(X, e_a)e_a, Y)= \sum_{a=1}^{2n}g(\check{R}(X, e_a)e_a, Y) 
\]
for all $X, Y\in\mc{D}$, where $\{e_a\}_{a=1}^{2n}$ is a $g$-orthonormal basis of $\mc{D}$ (and hence also $\check{g}$-orthonormal).  
As a consequence of  (\ref{transvR}), below  we derive  a higher-dimensional  generalization of an identity that satisfies the  transverse Ricci tensor on a Sasakian manifold (see part (ii) in \cite[Theorem 7.3.12]{BG08}, and be aware that here we mean the application of Theorem 7.3.12 to  the   Sasakian case).

Let us first derive a  simple formula for the Riemannian Ricci tensor, which will be useful at several points later in this work. So, let us fix  an $\mc{S}$-manifold $(M^{2n+s}, \phi, \xi_i, \eta_j, g)$.  As a  local $g$-orthonormal frame   $\{e_{a}\}_{a=1}^{2n+s}$ of $M^{2n+s}$ we may consider a local  $\phi$-adapted frame, i.e.,
 \[
 \{e_{a}\}_{a=1}^{2n+s}=\big\{E_1, \ldots, E_n, \phi(E_1), \ldots, \phi(E_n), \xi_1, \ldots, \xi_s\big\}
 \] 
with $\mc{D}={\rm span}\{E_k, \phi(E_k) : 1\leq k\leq n\}$ and $\mc{D}^{\perp}={\rm span}\{\xi_\ell : 1\leq\ell\leq s\}$, respectively.

\bl\label{Ricci_D}
Let $(M^{2n+s}, \phi, \xi_i, \eta_j, g)$ be an $\mc{S}$-manifold. 
Then, the Riemannian Ricci tensor satisfies 
 \begin{equation}\label{aster}
 \Ric^{g}(X, Y)=\Ric^{g}_{\mc{D}}(X, Y)+s g(X, Y)\,,
\end{equation}
 for any $X, Y\in\mc{D}$, where $\Ric^{g}_{\mc{D}}(X, Y)$ is defined by
 \begin{equation}\label{ricgD1}
 \Ric^{g}_{\mc{D}}(X, Y):=\sum_{k=1}^{n}\Big\{g(R^{g}(X, E_{k})E_{k}, Y)+g(R^{g}(X, \phi(E_k))\phi(E_k), Y)\Big\}\,.
 \end{equation}
\el
\pr
 By Proposition \ref{known1} it follows that   $ R^{g}(X, \xi_{\ell})\xi_{\ell}=-(\nabla^{g}_{X}\phi)\xi_{\ell}=\phi(\nabla^{g}_{X}\xi_{\ell})=\phi(-\phi(X))=X$, 
 for any $X\in\mc{D}$ and $\ell\in\{1, \ldots, s\}$ (alternatively, this  can be verified by a direct  computation  based on the identity $\nabla^{g}_{\xi_\ell}\phi=0$, see \cite{CFF90}). Thus, 
  \begin{eqnarray*}
 \Ric^{g}(X, Y)&=&\sum_{a=1}^{2n+s}R^{g}(X, e_a, e_a, Y)= \Ric^{g}_{\mc{D}}(X, Y)+\sum_{\ell=1}^{s}R^{g}(X, \xi_\ell, \xi_\ell, Y)\\
 &=& \Ric^{g}_{\mc{D}}(X, Y)+\sum_{\ell=1}^{s}g(R^{g}(X, \xi_\ell)\xi_\ell, Y)=\Ric^{g}_{\mc{D}}(X, Y)+\sum_{\ell=1}^{s}g(X, Y)\\
 &=&\Ric^{g}_{\mc{D}}(X, Y)+s g(X, Y)
 \end{eqnarray*}
 for any $X, Y\in\mc{D}$. This proves the given formula.
\pro

 \bc\label{transvRic}
 Let $(M^{2n+s}, \phi, \xi_i, \eta_j, g)$ be an $\mc{S}$-manifold. 
 The Ricci tensors $\Ric^g$ and $\Ric^{\check{g}}$ are related by
 \begin{equation}\label{Rictrans2}
 \Ric^{g}(X, Y)= \Ric^{\check g}(X, Y)-2s\,\check{g}(X, Y)\,,
 \end{equation}
 for all $ X, Y\in\mc{D}$. Therefore, 
 the scalar curvature $\Scal^g$ and the transverse scalar curvature $\Scal^{\check{g}}$ are such that
 \[
 \Scal^g=\Scal^{\check{g}}-2\,n\,s\,.
 \]
 In particular, for $s=1$ this gives the known identities  $\Ric^g|_{\mc{D}\times\mc{D}}=\Ric^{\check{g}}-2g|_{\mc{D}\times\mc{D}}$ and $\Scal^g=\Scal^{\check{g}}-2n$, respectively.
 \ec
 \pr
 We only present the proof of (\ref{Rictrans2}), using a local $g$-orthonormal $\phi$-adapted frame as above.  
 By the definition of the transverse Ricci tensor  $\Ric^{\check{g}}$ and  by (\ref{transvR}) 
we see that
\begin{eqnarray*}
\Ric^{g}_{\mc{D}}(X, Y)&=&\sum_{k=1}^{n}\Big\{g(R^{g}(X, E_{k})E_{k}, Y)+g(R^{g}(X, \phi(E_k))\phi(E_k), Y)\Big\}\\
&=&\Ric^{\check{g}}(X, Y)+s\sum_{k=1}^{n}\Big\{F(E_k, E_k)g(\phi(X), Y)-F(X, E_{k})g(\phi(E_k), Y)\\
&&+F(\phi(E_k), \phi(E_k))g(\phi(X), Y)-F(X, \phi(E_k))g(\phi^2(E_k), Y)\Big\}\\
&&-2s\sum_{k=1}^{n}\Big\{F(X, E_{k})g(\phi(E_k), Y)+F(X, \phi(E_k))g(\phi^2(E_k), Y)\Big\}\\
&=&\Ric^{\check{g}}(X, Y)-s\Big\{F(X, E_{k})g(\phi(E_k), Y)+F(X, \phi(E_k))g(\phi^2(E_k), Y)\Big\}\\
&&-2s\sum_{k=1}^{n}\Big\{g(X, \phi(E_{k}))g(\phi(E_k), Y)+g(X, E_k)g(E_k, Y)\Big\}\\
&=&\Ric^{\check{g}}(X, Y)-s\sum_{k=1}^{n}\Big\{g(X, \phi(E_k))g(Y, \phi(E_k))+g(X, E_k)g(Y, E_k)\Big\}\\
&&-2s\sum_{k=1}^{n}\Big\{g(X, \phi(E_{k}))g(Y, \phi(E_k))+g(X, E_k)g(Y, E_k)\Big\}\\
&=&\Ric^{\check{g}}(X, Y)-3s\sum_{k=1}^{n}\Big\{g(X, \phi(E_{k}))g(Y, \phi(E_k))+g(X, E_k)g(Y, E_k)\Big\}\\
&=&\Ric^{\check{g}}(X, Y)-3s\,g(X, Y)\,,
 \end{eqnarray*}
  for any $X, Y\in\mc{D}$. 
 The claim now follows by  (\ref{aster})  and the identification $g|_{\mc{D}\times\mc{D}}=\check{g}$.
 \pro

\subsection{The (generalized) $\eta$-Einstein condition and transverse K\"ahler-Einstein metrics}
Although  no $\mc{S}$-manifold $M^{2n+s}$ of CR-codimension $s>1$ can be Einstein, it is still meaningful to consider the notion of \textsf{$\eta$-Einstein $\mc{S}$-manifolds}, thereby extending the notion of an $\eta$-Einstein Sasakian manifold.  We   recall the following 
\bd\textnormal{(\cite[pp.~432-433]{KT72})}\label{etaE}
An $\mc{S}$-manifold $(M^{2n+s}, \phi, \xi_i, \eta_j, g)$   is said to be $\eta$-Einstein when the Riemannian Ricci tensor $\Ric^g$ satisfies the condition 
\begin{equation}\label{etaE1}
\Ric^{g}(X, Y)=\al\, g(X, Y)+\beta \sum_{i=1}^{s}\eta_{i}(X)\eta_{i}(Y)+(\al+\beta)\sum_{1\leq i\neq j\leq s}\eta_{i}(X)\eta_{j}(Y)
\end{equation}
for some constants $\al, \beta$ and all $X, Y\in TM$.
\ed
In compact form, the $\eta$-Einstein condition may be expressed as  
\[
\Ric^{g}=\al\cdot g+\beta\cdot\sum_{i}\eta_{i}\otimes\eta_{i}+(\al+\beta)\cdot\sum_{i\neq j}\eta_{i}\otimes\eta_{j}\,.
\]
 Obviously, for $s=1$ this defines an $\eta$-Einstein Sasakian manifold, see   \cite[Chapter 11]{BG08} for details. 
 Based on (\ref{not_diag}) and (\ref{etaE1})  one can  prove that the constants $\al, \beta$ are such that $\al+\beta=2n$.\footnote{Note that a different normalization is used for relation (\ref{not_diag}) in \cite{KT72}.}   Concerning the scalar curvature $\Scal^g$ of an $\eta$-Einstein $\mc{S}$-manifold,  one has 
\begin{eqnarray*}
\Scal^g&=&\sum_{a=1}^{2n+s}\Ric^{g}(e_a, e_a)=\al\sum_{a=1}^{2n+s}g(e_a, e_a)+\beta\sum_{k=1}^{s}\sum_{i=1}^{s}\eta_{i}(\xi_k)\eta_{i}(\xi_k)+(\al+\beta)\sum_{k=1}^{s}\sum_{i\neq j}\eta_{i}(\xi_k)\eta_{j}(\xi_k)\\
&=&\al\cdot (2n+s)+\beta\cdot s+0=2\cdot\al\cdot n+s\cdot (\al+\beta)=2n\cdot (\al+s)\,,
\end{eqnarray*}
for any $s\geq 1$. Hence
\bc
 Any $\eta$-Einstein $\mc{S}$-manifold is of constant scalar curvature  $\Scal^{g}=2n(\al+s)$, where $\al$ is the constant determined by the $\eta$-Einstein condition {\rm(\ref{etaE1})}.
 \ec
For an $\eta$-Einstein Sasakian manifold $(M^{2n+1}, \phi, \xi, \eta, g)$ this gives $\Scal^g=2n(\al+1)$, hence it recovers Corollary 9 in \cite{BGM96}.
Let us finally extend a familiar definition from the Sasakian case to the setting of $\mc{S}$-manifolds (cf. \cite[Definition 3.3]{FOW09} for $s=1$).
\bd\label{transvKE}
An $\mc{S}$-manifold $(M^{2n+s}, \phi, \xi_i, \eta_j, g)$ is said to be \textsf{transverse K\"ahler-Einstein}
if the transverse Ricci tensor $\Ric^{\check{g}}$ is a multiple of the transverse metric $\check{g}$, i.e., 
$\Ric^{\check{g}}=\lambda\cdot\check{g}$ for some real constant $\lambda$.
\ed
We will now prove the main result of this section which provides a natural bijection between the  notion of  $\eta$-Einstein $\mc{S}$-manifolds and the transverse K\"ahler-Einstein condition.
\bt\label{Main_Corres}
Let $(M^{2n+s}, \phi, \xi_i, \eta_j, g)$ be an $\eta$-Einstein $\mc{S}$-manifold. Then,  $M$ is transverse K\"ahler-Einstein with Einstein constant $\lambda=\al+2s$. Conversely, if $M$ is transverse K\"ahler-Einstein with Einstein constant $\lambda$, then $M$ is $\eta$-Einstein with constants $\al, \beta$ given by $\al=\lambda-2s$ and $\beta=2n-(\lambda-2s)$, respectively. \et
\pr
 Suppose that     $\Ric^{g}=\al\cdot g+\beta\cdot \sum_{j=1}^{s}\eta_{j}\otimes\eta_{j}+(\al+\beta)\cdot\sum_{i\neq j}\eta_{i}\otimes\eta_{j}$ 
    for some constants $\al, \beta$, with $\al+\beta=2n$. For $X, Y\in\mc{D}$ we   obtain that 
    \begin{equation}\label{Ricgal}
    \Ric^{g}(X, Y)=\al\,g(X, Y)=\al\,g(\phi X, \phi Y)=\al\,\check{g}(X, Y)\,,
    \end{equation}
     where $\check{g}$ is the transverse K\"ahler metric defined by (\ref{TKahlerg}). Then, the relation 
     \[
     \Ric^{g}(X, Y)=\Ric^{\check{g}}(X, Y)-2\,s\,\check{g}(X, Y)\,,
     \]
      for all $X, Y\in\mc{D}$ proved in Corollary \ref{transvRic},   reduces to $\Ric^{\check{g}}=(\al+2\,s)\check{g}$.  This proves the one direction.\\
      \noindent
      Conversely, suppose that $\Ric^{\check{g}}=\lambda\cdot\check{g}$ for some real constant $\lambda$. Then, again by   Corollary \ref{transvRic}, it follows that $\Ric^{g}(X, Y)=(\lambda-2s)g(X, Y)$, 
 for all $X, Y\in\mc{D}$.
 Since by (\ref{not_diag}) we should also have $\Ric^{g}(\xi_k, \xi_\ell)=2n$ for all $k, \ell\in\{1, \ldots, s\}$,  it follows that we may expresses $\Ric^g$ as 
 \[
 \Ric^{g}=(\lambda-2s)g+\beta\sum_{i}\eta_{i}\otimes\eta_{i}+2n\sum_{i\neq j}\eta_{i}\otimes\eta_{j}\,,
 \]
with $\beta=2n-(\lambda-2s)$. This proves our claim.
\pro
In this way, we establish a bijective correspondence between the $\eta$-Einstein condition and the transverse K\"ahler-Einstein condition, thereby generalizing  the  well-known  situation for $s=1$, i.e., the Sasakian case (cf. \cite[Theorem 11.1.3]{BG08} or \cite{FOW09}). 

 \br
{\rm (a)} 
 Let $M^{2n+s}$ be an $\eta$-Einstein $\mc{S}$-manifold.  By Theorem \ref{Main_Corres} we have the relation $\Ric^{\check{g}}=(\al+2\,s)\check{g}$. Hence   $\Ric^{\check{g}}$ is positive-definite if and only if $\al+2s>0$. If $\al+2s=0$, then the transverse metric is $\Ric^{\check{g}}$-flat.  \\
{\rm(b)} Let $\rho(X, Y):=\Ric^{\check{g}}(X, \check{J}Y)$ be the \textsf{transverse Ricci 2-form} induced by $\Ric^{\check{g}}$. Then the condition $\Ric^{\check{g}}=\lambda\check{g}$ gives $\rho=\lambda\check\om$. Thus, when $M^{2n+s}$ is an $\eta$-Einstein $\mc{S}$-manifold it follows that  
$\rho=(\al+2\,s)\hat{\om}=(\al+2\,s)F$. 
For $s=1$ the transverse Ricci form was used in \cite{BGM96}  to define the \textsf{basic first Chern class} of $\mc{D}$. We will address the extension of this construction in our higher-dimensional scenario separately, in a forthcoming work, together with applications of the so-called (generalized) \textsf{$\mc{D}$-homothetic transformations} (see \cite{CDT07}).
 \er

 
\section{Curvature of the characteristic connection on $\mc{S}$-manifolds}\label{CSmnfds}

\subsection{Metric connections with parallel skew-torsion}\label{MCST}
Let $(M, g)$ be a Riemannian manifold and let $\nabla^g$ be the Levi-Civita connection associated to $g$. 
 Recall that the torsion $T$ of a metric connection $\nabla$ on $M$  is the vector-valued 2-form defined by
 $T(X, Y)=\nabla_{X}Y-\nabla_{Y}X-[X, Y]$, for all $X, Y\in TM$. 
The connection $\nabla$     is said to have \textsf{totally skew-symmetric torsion}, or \textsf{skew-torsion} in short,
whenever 
  \[
 T(X, Y, Z):=g(T(X, Y), Z)=-g(T(X, Z), Y)=-T(X, Z, Y)\,,
  \]
for all $X, Y, Z\in TM$. This means that the  induced tensor field  via $T$ by a contraction with $g$ is a 3-form on $M$.
In this case  $\nabla$ is expressed  as $\nabla=\nabla^{g}+\frac{1}{2}T$, that is, 
\[
g(\nabla_{X}Y, Z)=g(\nabla^{g}_{X}Y, Z)+\frac{1}{2}g(T(X, Y), Z)=g(\nabla^{g}_{X}Y, Z)+\frac{1}{2}T(X, Y, Z)\,,\quad X, Y, Z\in TM\,.
\]  
\br  
For any 3-form $T$ on $(M^n, g)$ we denote by   $\Ker(T)=\{X\in TM : X\lrcorner T=0\}\subset TM$ its kernel. 
Moreover,  we  assign to $T$   a  4-form $\sigma_T$ given by (cf. \cite{FrIv})
\[
\sigma_{T}:=\frac{1}{2}\sum_{i=1}^{n}(e_i\lrcorner T)\wedge(e_i\lrcorner T)\,,
\]
where $e_1, \ldots, e_n$ is a  local orthonormal frame of  $(M, g)$. This  is also expressed by
 \[
 \sigma_{T}(X, Y, Z, W)=\fr{S}_{X, Y, Z}T(T(X, Y), Z, W)=\fr{S}_{X, Y, Z}g(T(X, Y), T(Z, W))\,, 
 \]
for any $X, Y, Z, W\in TM$, where $\fr{S}_{X, Y, Z}$ denotes the cyclic sum over the vector fields $X, Y, Z\in TM$.   
For the curvature tensor $R^{\nabla}$ induced by $\nabla$  we adopt the sign convention
\[
R^{\nabla}(X, Y)Z=\nabla_{X}\nabla_{Y}Z-\nabla_{Y}\nabla_{X}Z-\nabla_{[X, Y]}Z\,,\quad X, Y, Z\in TM\,.
\]
\er

In this text we are mainly  interested in metric connections with parallel skew-torsion, 
\[
\nabla=\nabla^{g}+\frac{1}{2}T\,,\qquad \nabla T=0\,.
\]
Under the assumption $\nabla T=0$ we have  $\dd T=2\sigma_T$ and  the first  Bianchi identity for $R^{\nabla}$ becomes (cf. \cite{IP01, FrIv})
\[
\fr{S}_{X, Y, Z}R^{\nabla}(X, Y, Z, W)=\sigma_{T}(X, Y, Z, W)=\frac{1}{2}\dd T(X, Y, Z, W)\,,\quad X, Y, Z, W\in TM\,.
\]
This implies that the 4-tensor $R^{\nabla}(X, Y, Z, W)=g(R^{\nabla}(X, Y)Z, W)$ is symmetric  after interchanging  the blocks $(X, Y)$ and $(Z, W)$, i.e.,
\[
R^{\nabla}(X, Y, Z, W)=R^{\nabla}(Z, W, X, Y)\,,\quad X, Y, Z, W\in TM\,.
\]
Recall that the curvature tensor $R^{g}$ associated to   $\nabla^g$  satisfies the same symmetry property.

Next we will encode the above situation by tripes $(M^{n}, g, T)$ and assume, once and for all,  that $M$ is connected and oriented. Since we are interested in the case where $\nabla T=0$, we will henceforth refer to  $(M^{n}, g, T)$ as a \textsf{geometry with parallel skew-torsion}. As usual, we will identify the tangent space $(T_xM, g_x)$ of $M$ at an arbitrary point $x\in M$ with the Euclidean space  $(\R^n, ( \ , \ ))$. The  holonomy group of $\nabla:=\nabla^{g}+\frac{1}{2}T$, denoted by $\Hol(\nabla)$, can be then viewed as  a subgroup of $\SO(n)$, and its Lie algebra (holonomy algebra) $\fr{hol}(\nabla)$ is a  Lie subalgebra of $\fr{so}(n)=\Lie(\SO(n))$.  Recall that there is a natural action of $\Hol(\nabla)$ on the space of tensors fields and in particular on $TM$, referred to as the holonomy representation.   When $(M^{n}, g, T)$ admits a $\Gg$-structure for some closed subgroup $\Gg\subset\SO(n)$ and $\nabla$ is an adapted connection for this $\Gg$-structure, then, $\nabla$ is referred to as the \textsf{characteristic connection}.
In this case, 
as a consequence of the general theory of $\Gg$-structures, we have the inclusion $\Hol(\nabla)\subset \Gg\subset\SO(n)$.  Next, by slight abuse of notation, given any geometry $(M^n, g, T)$ with parallel skew-torsion  we will  always call the connection $\nabla=\nabla^{g}+\frac{1}{2}T$ characteristic.

Another group having a distinguished role for a geometry with parallel skew-torsion    is the \textsf{stabilizer} $\Stab(T)=\{A\in\SO(n) : A\cdot T=T\}$ 
of the torsion 3-form  $T$ in $\SO(n)$.  The geometric significance of $\Stab(T)$ follows from the fact that when  the condition  $\nabla T=0$ is satisfied, then the holonomy group $\Hol(\nabla)$ should preserve $T$ and thus is a subgroup of the stabilizer of $T$, i.e.,  $\Hol(\nabla)\subset\Stab(T)$.  Hence the Lie algebra $\fr{stab}(T)$ of $\Stab(T)$  is  a Lie subalgebra of $\fr{so}(n)$ such that  $\fr{hol}(\nabla)\subset\fr{stab}(T)\subset\fr{so}(n)$. The stabalizer of $T$ has played a role in many works, see for example \cite{Schoe, AFF15, MS24}.


 \subsection{Adapted connections with skew-torsion on metric $f$-manifolds}
Let us now consider metric $f$-manifolds whose characteristic vector fields commute, and recall the following 
result proved in \cite{BC} (for the special cases  $s=0, 1$, Theorem \ref{ThmBCH} was known by   \cite{FrIv}).
\bt\textnormal{(\cite{BC})}\label{ThmBCH}
Let  $(M^{2n+s},  \phi, \xi_i, \eta_j, g)$  be a metric $f$-manifold  
 whose characteristic
vector fields  $\xi_i$ commute, i.e., $[\xi_i, \xi_j]=0$ for all $i, j\in\{1, \ldots, s\}$. Then the following conditions are equivalent: \\
(1) $N^{(1)}$ is totally skew-symmetric and  $\xi_1, \ldots, \xi_s$ are Killing   vector fields.\\ 
(2)    $(M^{2n+s},  \phi, \xi_i, \eta_j, g)$  admits a  metric  connection with skew-torsion 
$T$ preserving the structure, i.e.,
\begin{equation}\label{parallelall}
\nabla\phi=0\,,\quad \nabla \xi_i=0\,,\quad \nabla\eta_j=0\,,\quad\text{for all} \quad i, j\in\{1, \ldots, s\}\,.
\end{equation}
In this case the connection $\nabla$ is uniquely determined by the relation $\nabla=\nabla^g+\frac{1}{2}T$, where the torsion 3-form $T$   is   given by  
\begin{equation}\label{TorsionBC}
T=\sum_{i=1}^{s}\eta_{i}\wedge\dd\eta_i+\dd^{\phi}F+N^{(1)}-\sum_{i=1}^{s}(\eta_{i}\wedge(\xi_i\lrcorner N^{(1)}))\,,
\end{equation}
with $\dd^{\phi}F:=-\dd F\circ\phi$.
\et

In the present article we will mainly restrict our attention   to $\mc{S}$-manifolds $(M^{2n+s},  \phi, \xi_i, \eta_j, g)$.   
\bt \textnormal{(\cite{BC})}\label{Torsion_skewS}
A  contact metric $f$-manifold $(M^{2n+s},  \phi, \xi_i, \eta_j, g)$ admits  a metric  connection $\nabla$ with skew-torsion  satisfying the conditions in  {\rm(\ref{parallelall})} if and only if $M^{2n+s}$ is an $\mc{S}$-manifold. 
The connection $\nabla$  is then uniquely determined by the relation $\nabla=\nabla^g+\frac{1}{2}T$, where the torsion 3-form $T$ given in {\rm(\ref{TorsionBC})} reduces to
\[
T=\sum_{i=1}^{s}\eta_{i}\wedge\dd\eta_i\,.
\]
Moreover, in this case  the torsion 3-form $T$   is $\nabla$-parallel, i.e., $\nabla T=0$.
\et
For $s=1$ the corresponding characteristic connection and the condition  $\nabla T=0$ was first presented in   \cite{FrIv}.  In \cite{BC} we also proved that 
\bp\textnormal{(\cite{BC})}
 On an $\mc{S}$-manifold    $(M^{2n+s},  \phi, \xi_i, \eta_j, g)$ the kernel of the torsion 3-form $T$   is a smooth  subbundle of the vertical distribution $\mc{D}^{\perp}$ of rank $s-1$. 
\ep
Hence only for Sasakian manifolds we have $\Ker(T)=\{0\}$.
Let us also recall the expression of  the vector-valued   2-form $T$ on an $\mc{S}$-manifold (here at some place we combine \cite[Proposition 4.4]{BC} with Lemma \ref{XYvert}).

\bp\textnormal{(\cite{BC})}\label{2-valuedT}
Let $(M^{2n+s}, \phi,\xi_i, \eta_j, g)$ be an $\mc{S}$-manifold endowed with its  characteristic connection $\nabla$ with skew-torsion $T$, as presented in Theorem \ref{Torsion_skewS}.   Then  
\[
T(X, Y)=2\sum_{i=1}^{s}\big\{F(X, Y)\xi_i-\eta_{i}(X)\phi(Y)+\eta_i(Y)\phi(X)\big\}
\]
 for all $X, Y\in TM$. Thus,  
 \[
T(X, \xi_i)=2\phi(X)\,, \quad T(\xi_i, \xi_j)=0\,,\quad T(X, Y)=2\sum_{i=1}^{s}F(X, Y)\xi_{i} =2F(X, Y)\bar{\xi}=-[X, Y]_{\mc{D}^{\perp}}\,,
\]
for all $X, Y\in\mc{D}$ and $i, j\in\{1, \ldots, s\}$, where $\bar{\xi}:=\sum_{i=1}^{s}\xi_i$.  It also follows that  $T(X, Y, Z)=0$ for all $X, Y, Z\in\mc{D}$.
\ep

Finally, let us mention that since $\phi|_{\mc{D}}$ is an almost complex structure on the horizontal part $\mc{D}$,  and the structure group of an $\mc{S}$-manifold $(M^{2n+s}, \phi, \xi_i, \eta_j, g)$ is the Lie group $\U(n)\times\Id_{s}\subset\SO(2n+s)$, it follows that
 \bp\label{stabilizerS}
If   $(M^{2n+s}, \phi,\xi_i, \eta_j, g)$  is an $\mc{S}$-manifold, then the torsion 3-form $T$ of its characteristic connection $\nabla$  has stabilizer  $\fr{stab}(T)$    in $\fr{so}(2n+s)$  isomorphic to $\fr{u}(n)=\Lie(\U(n))$, i.e., $\fr{stab}(T)\cong\fr{u}(n)$.  Hence  $\fr{hol}(\nabla)\subset\fr{u}(n)$. 
\ep 
Note that  for the Sasakian case ($s=1$), this  is a  standard fact (see \cite[Section 2.11]{MS24}).


\subsection{Curvature properties}
 Let  $(M^{2n+s}, \phi,\xi_i, \eta_j, g)$  be an $\mc{S}$-manifold endowed with its  characteristic connection $\nabla$ described above. We will derive the curvature tensor $R^{\nabla}$ of $\nabla$.

\bp\label{pc1}
On an $\mc{S}$-manifold  $(M^{2n+s}, \phi,\xi_i, \eta_j, g)$   the curvature tensor $R^{\nabla}$  of the characteristic connection $\nabla$  satisfies  the relation $R^{\nabla}(X, Y, Z, W)=R^{g}(X, Y, Z, W)+\Sigma_{s}(X, Y, Z, W)$,  for all $X, Y, Z, W\in\mc{D}$, where 
 \[
 \Sigma_{s}(X, Y, Z, W):=2s\Big\{F(X, Y)F(Z, W)+\frac{1}{2}F(Y, Z)F(X, W)+\frac{1}{2}F(Z, X)F(Y, W)\Big\}\,.
 \]
All the other combinations vanish,  and in particular we have that
\begin{align*}
R^{\nabla}(X, Y, \xi_i, Z)&=R^{\nabla}(\xi_i, Z, X, Y)=-R^{\nabla}(Z, \xi_i, X, Y) =-R^{\nabla}(X, Y, Z, \xi_i)=0\,, \\
R^{\nabla}(X, \xi_i, Y, \xi_j)&=R^{\nabla}(X, Y, \xi_i, \xi_j)=R^{\nabla}(\xi_i, \xi_j, X, Y)=0\,,\\
R^{\nabla}(X, \xi_i, \xi_j, \xi_k)&=0\,,\quad R^{\nabla}(\xi_i, \xi_j, \xi_k, \xi_\ell)=0\,,
\end{align*}
for all $ X, Y, Z\in\mc{D}$ and  $i, j, k, \ell\in\{1, \ldots, s\}$
\ep
\pr
 Since $\nabla T=0$,  by formula (3.19) in \cite{IP01} and the definition of $\sigma_T$, it follows that  the curvature tensor $R^{\nabla}$ of $\nabla$ is such that
 \begin{equation}\label{c1}
R^{\nabla}(X, Y, Z, W)=R^{g}(X, Y, Z, W)+\frac{1}{4}g(T(X, Y), T(Z, W))+\frac{1}{4}\sigma_{T}(X, Y, Z, W)\,,
\end{equation}
for all $X, Y, Z, W\in  TM$. In addition,  $R^{\nabla}$ has the usual symmetry properties as $R^{g}$, which means that: 
$ (a) \ R^{\nabla}(X, Y, Z, W)=-R^{\nabla}(Y, X, Z, W)$, $(b) \ R^{\nabla}(X, Y, Z, W)=-R^{\nabla}(X, Y, W, Z)$ and $(c) \ R^{\nabla}(X, Y, Z, W)=R^{\nabla}(Z, W, X, Y)$, 
for all $X, Y, Z, W\in TM$. Note that $(a)$ and $(b)$ hold for any metric connection, while the third relation (block symmetry) is  a consequence of the parallelism of $T$, as   mentioned above.

\noindent For $X, Y, Z, W\in\mc{D}$ by  Proposition  \ref{2-valuedT} we have $T(X, Y)=2\sum_{i}F(X, Y)\xi_i$ and $T(Z, W)=2\sum_{j}F(Z, W)\xi_j$,  
and thus in this case  $\sigma_T$ does not vanish, i.e.,
\[
\sigma_{T}(X, Y, Z, W)=4s\left\{F(X, Y)F(Z, W)+F(Y, Z)F(X, W)+F(Z, X)F(Y, W)\right\}\,.
\]
The given relation for $R^{\nabla}(X, Y, Z, W)$ now follows  by (\ref{c1}).\\
 \noindent   Let $X, Y, Z\in\mc{D}$ and $W=\xi_i\in\mc{D}^{\perp}$. In this case it is easy to see that $ \sigma_{T}(X, Y, Z,  \xi_i)=0$.  In fact,  the distributions $\mc{D}$ and $\mc{D}^{\perp}$ are $\nabla$-parallel, hence we  have $g(R^{\nabla}(X, Y)\mc{D}, \mc{D}^{\perp})=0$   for all  $X, Y\in TM$.  Therefore $R^{\nabla}(X, Y, Z, \xi_i)=0$ and by  the symmetry mentioned in part (b)  above, we also have $R^{\nabla}(X, Y, \xi_i, Z)=0$.  An alternative is based on Proposition \ref{known1}  and the relation (\ref{c1}), which we avoid to present.\\
 \noindent Let $X, Y\in\mc{D}$. Since $T(\xi_i, \xi_j)=0$ for all $i, j\in\{1, \ldots, s\}$ and $T(X, \xi_i)=2\phi(X)$ for all $i\in\{1, \ldots, s\}$,  we see that
 \begin{eqnarray*}
 \sigma_{T}(X, \xi_i, Y, \xi_j)&=&g(T(X, \xi_i), T(Y, \xi_j))+g(T(\xi_i, Y), T(X, \xi_j))+g(T(Y, X), T(\xi_i,  \xi_j))\\
 &=&4g(\phi(X), \phi(Y))-4g(\phi(Y), \phi(X))=0\,.
 \end{eqnarray*}
 Thus, by Proposition   \ref{known1}, the relation (\ref{c1}) becomes
 \begin{eqnarray*}
 R^{\nabla}(X, \xi_i, Y, \xi_j)&=&g(R^{g}(X, \xi_i)Y, \xi_j)+\frac{1}{4}g(T(X, \xi_i), T(Y, \xi_j)\\
 &=&-g((\nabla^g_{X}\phi)Y, \xi_j)+g(\phi(X), \phi(Y))\,.\quad (\ddag)
 \end{eqnarray*}
Now, for the first term in this expression we may use  relation (2.6)  in \cite{BC}, which gives
 \begin{eqnarray*}
 -g((\nabla^g_{X}\phi)Y, \xi_j)&=&g((\nabla^{g}_{X}\phi)\xi_j, Y)=-g(\phi(\nabla^{g}_{X}\xi_j), Y)\\
 &\overset{\rm(\ref{gphi})}{=}& g(\nabla^{g}_{X}\xi_j, \phi(Y))=-g(\phi(X), \phi(Y))
 \end{eqnarray*}
 since $(\nabla^{g}_{X}\phi)\xi_j=-\phi(\nabla^{g}_{X}\xi_j)$ and  moreover on an  $\mc{S}$-manifold we have $\nabla^{g}_{X}\xi_j=-\phi(X)$ for all $X\in TM$ and $j\in\{1, \ldots, s\}$.   Thus by $(\ddag)$ we get  $ R^{\nabla}(X, \xi_i, Y, \xi_j)=0$ for all $X, Y\in\mc{D}$ and $i, j\in\{1, \ldots, s\}$.  Similarly one can prove that  $R^{\nabla}(X, Y, \xi_i, \xi_j)=0$ for $X, Y\in\mc{D}$ and $i, j\in\{1, \ldots, s\}$, and the rest relations follow by the symmetry properties $(a)$-$(c)$ mentioned in the beginning of the proof.\\
 \noindent  Let $X\in\mc{D}$.  Since  $\nabla^{g}_{\xi_i}\xi_j=0=[\xi_i, \xi_j]$ for all $i, j\in\{1, \ldots, s\}$ it is easy to see that $R^g(\xi_j, \xi_k)\xi_i=0$ and it  follows that  $R^{g}(\xi_j, \xi_k,\xi_i, X)=- R^{g}(\xi_j, \xi_k, X, \xi_i)=-R^{g}(X, \xi_i, \xi_j, \xi_k)=0$, 
for all $X\in\mc{D}$. Then, since $T(\xi_i, \xi_j)=0$ for all $i, j\in\{1, \ldots, s\}$   it is also easy to see that $R^{\nabla}(X, \xi_i, \xi_j, \xi_k)=0$. The final case is treated similarly. 
\pro

\bc\label{curvatCorl}
On an $\mc{S}$-manifold $(M^{2n+s}, \phi,\xi_i, \eta_j, g)$ the curvature tensor $R^{\nabla}$  of the characteristic connection $\nabla$  satisfies
\begin{equation}\label{cc1}
R^{\nabla}(X, Y)Z=R^{g}(X, Y)Z-2sF(X, Y)\phi(Z)-sF(Y, Z)\phi(X)-sF(Z, X)\phi(Y)\,,
\end{equation}
for any $X, Y, Z\in\mc{D}$.  Moreover, 
\begin{equation}\label{useRic1}
R^{\nabla}(\xi_i, \xi_j)\xi_k= R^{\nabla}(X, \xi_i)\xi_j = R^{\nabla}(X, Y)\xi_i =
R^{\nabla}(X, \xi_i)Y = 
R^{\nabla}(\xi_i, \xi_j)X=0 
 \end{equation}
for all $X, Y\in\mc{D}$ and $i, j, k\in\{1, \ldots, s\}$
\ec


 \subsection{The Ricci tensor and the scalar curvature}\label{charRicScal}
We will now present the Ricci tensor $\Ric^{\nabla}$ induced by $\nabla$.
Based on the condition  $\nabla T=0$ one can show that $\delta^gT=0$ and hence  the Ricci tensor $\Ric^{\nabla}$ associated to $\nabla=\nabla^g+\frac{1}{2}T$  is   symmetric. This is defined by  
\begin{equation}\label{ricci1}
\Ric^\nabla(U,V)=\sum_{a=1}^{2n+s}g\big(R^{\nabla}(U, e_a)e_a, V\big)\,,\quad U,V \in TM\,,
\end{equation}
where $\{e_a\}_{a=1}^{2n+s}$ is a local $g$-orthonormal frame of $(M^{2n+s}, \phi,\xi, \eta, g)$.
As in the previous section,  one may assume that $\{e_{a}\}$ is a local  $\phi$-basis, i.e.,
\[
\{e_{a}\}_{a=1}^{2n+s}=\big\{E_1, \ldots, E_n, \phi(E_1), \ldots, \phi(E_n), \xi_1, \ldots, \xi_s\big\}\,.
\]
 %

\bp\label{ricp1}
Let   $(M^{2n+s}, \phi,\xi_i, \eta_j, g)$ be an $\mc{S}$-manifold. Then the Ricci tensor $\Ric^{\nabla}$  induced by the characteristic connection $\nabla=\nabla^g+\frac{1}{2}T$  satisfies the following relations:
\begin{eqnarray*}
\Ric^{\nabla}(X, Y)&=&\Ric^{g}(X, Y)-2\,s\,g(X, Y)\,,\quad \forall \ X, Y\in\mc{D}\,,\\
\Ric^{\nabla}(X, \xi_i)&=&0\,,\quad  \forall \ X\in\mc{D},\ i\in\{1, \ldots, s\}\,,\\
\Ric^{\nabla}(\xi_i, \xi_j)&=&0\,,\quad\forall \ i, j\in\{1, \ldots, s\}\,.
\end{eqnarray*}
In particular,  for the Ricci endomorphism $\Ric^{\nabla} : TM\to TM$ corresponding to $\nabla$  we obtain 
\[
\Ric^{\nabla}|_{\mc{D}}=\Ric^{g}|_{\mc{D}}-2\,s\Id|_{\mc{D}}\,,\qquad \Ric^{\nabla}|_{\mc{D}^{\perp}}=0\,,
\]
where $\Ric^{g}|_{\mc{D}}$  denotes the restriction of the Riemannian Ricci endomorphism to the horizontal distribution $\mc{D}$.
\ep

\pr
By (\ref{ricci1}) and   Proposition \ref{pc1}  we see that
\[
\Ric^{\nabla}(X, \xi_i)=\sum_{k=1}^{n}\{R^{\nabla}(X, E_k,   E_k, \xi_i)+R^{\nabla}(X, \phi(E_k), \phi(E_k), \xi_i)\}+\sum_{\ell=1}^{s}R^{\nabla}(X, \xi_\ell,  \xi_\ell, \xi_i)=0\,.
\]
 Similarly, we get $\Ric^{\nabla}(\xi_i, \xi_j)=0$ for all $i, j\in\{1, \ldots, j\}$. Suppose now that $X, Y$ are two horizontal vector fields, i.e., $X, Y\in\mc{D}$.
 Then
 \begin{eqnarray*}
 \Ric^{\nabla}(X, Y)&=&\sum_{k=1}^{n}\{R^{\nabla}(X, E_k, E_k, Y)+R^{\nabla}(X, \phi(E_k),   \phi(E_k), Y)\}+\sum_{\ell=1}^{s}R^{\nabla}(X, \xi_\ell,  \xi_\ell, Y)\\
 &=&\sum_{k=1}^{n}g(R^{\nabla}(X, E_{k})E_{k}, Y)+\sum_{k=1}^{n}g(R^{\nabla}(X, \phi(E_k))\phi(E_k), Y)
 \end{eqnarray*}
 since $R^{\nabla}(X, \xi_i, \xi_j, Y)=0$ for all $X, Y\in\mc{D}$ and $i, j\in\{1, \ldots, s\}$. Now, by (\ref{cc1}) and since $F$ is a 2-form, a direct computation shows that
 \begin{eqnarray*}
 R^{\nabla}(X, E_k)E_k&=&R^{g}(X, E_k)E_k-sF(X, E_k)\phi(E_k)\,,\\
 R^{\nabla}(X, \phi E_k)\phi E_k&=&R^{g}(X, \phi E_k)\phi E_k +sF(X, \phi E_k)E_{k}\,.
 \end{eqnarray*}
 Thus, taking the sums we see that
 \begin{eqnarray*}
  \Ric^{\nabla}(X, Y)&=&\sum_{k=1}^{n}g(R^{\nabla}(X, E_{k})E_{k}, Y)+\sum_{k=1}^{n}g(R^{\nabla}(X, \phi(E_k))\phi(E_k), Y)\\
  &=&\sum_{k=1}^{n}\Big\{g(R^{g}(X, E_{k})E_{k}, Y)+g(R^{g}(X, \phi(E_k))\phi(E_k), Y)\Big\}\\
  &&-s\sum_{k=1}^{n}\Big\{F(X, E_{k})g(\phi(E_k), Y)-F(X, \phi(E_k))g(E_{k}, Y)\Big\}\\
  &=&\Ric^{g}_{\mc{D}}(X, Y)-s\sum_{k=1}^{n}\Big\{F(X, E_{k})g(\phi(E_k), Y)-F(X, \phi(E_k))g(E_{k}, Y)\Big\}\\
  &=&\Ric^{g}_{\mc{D}}(X, Y)-s\sum_{k=1}^{n}\Big\{g(X, \phi(E_k))g(Y, \phi(E_k))+g(X, E_{k})g(Y, E_k)\Big\}\\
 &=&\Ric^{g}_{\mc{D}}(X, Y)-sg(X, Y)\,,  \quad\quad\quad\quad\quad\quad\quad\quad\quad\quad\quad\quad\quad\quad\quad\quad\quad\quad\quad\quad (\ast)
 \end{eqnarray*}
 for any $X, Y\in\mc{D}$, where $ \Ric^{g}_{\mc{D}}(X, Y)$ is defined by (\ref{ricgD1}).
 Now recall by (\ref{aster}) that  $\Ric^{g}_{\mc{D}}(X, Y)=\Ric^{g}(X, Y)-sg(X, Y)$, 
 for any $X, Y\in\mc{D}$,   and our claim follows by $(\ast)$.\\
 \pro

The component $\Ric^{g}_{\mc{D}}$ of the Riemannian Ricci tensor given in (\ref{ricgD1}) contains the   transverse Ricci tensor associated to the transverse metric $\check{g}$  of the characteristic foliation $\mc{F}$ defined on  an $\mc{S}$-manifold $M^{2n+s}$.  The precise relationship was established in  the proof of  Corollary \ref{transvRic}. Therefore, on the horizontal distribution $\mc{D}$,   both  Ricci tensors $\Ric^{g}$ and  $\Ric^{\nabla}$  may be described in terms of the transverse Ricci tensor $\Ric^{\check{g}}$ and  the transverse metric $\check{g}$, see for example the relation (\ref{chartransverse}) given below.   An analogous observation applies to the scalar curvatures  $\Scal^g$ and $\Scal^{\nabla}$, which are determined by the transverse scalar curvature $\Scal^{\check{g}}$ and the parameters $n, s$.  Below,  the relationship between the scalar curvatures $\Scal^g$ and $\Scal^{\nabla}$ is specified.

 \bc\label{scalar_curvature}
 On an $\mc{S}$-manifold  $(M^{2n+s}, \phi,\xi_i, \eta_j, g)$   the scalar curvature $\Scal^{\nabla}$ of the characteristic connection $\nabla$
 is related to the  Riemannian scalar curvature $\Scal^g$ by the formula 
 \[
 \Scal^{\nabla}=\Scal^g-6\,n\,s\,.
 \]
 \ec
  \pr
 The Ricci tensor $\Ric^\nabla$ is symmetric and  by definition
\[
\Scal^\nabla
=\sum_{a=1}^{2n+s}\Ric^{\nabla}(e_a, e_a)=\sum_{a=1}^{2n}\Ric^\nabla(e_a,e_a)+\sum_{\ell=1}^{s}\Ric^\nabla(\xi_\ell,\xi_\ell)=\sum_{a=1}^{2n}\Ric^\nabla(e_a,e_a)
\]
since   $\Ric^{\nabla}(\xi_\ell, \xi_\ell)=0$. 
Using the identity  $\Ric^\nabla(X,Y)=\Ric^g(X,Y)-2\,s\,g(X,Y)$ for any  $X,Y\in\mc{D}$ proved above, it follows that
 \begin{eqnarray*}
\Scal^\nabla&=&\sum_{a=1}^{2n}\Ric^\nabla(e_a,e_a)=\sum_{k=1}^{n}\Big\{\Ric^{\nabla}(E_k,E_k)+\Ric^{\nabla}(\phi(E_k),\phi(E_k))
\Big\} \\
&=&
\sum_{k=1}^{n}\Big\{\Ric^{g}(E_k,E_k)+\Ric^{g}(\phi(E_k),\phi(E_k))\Big\} - 2s\sum_{k=1}^{n}\Big\{g(E_k,E_k)+g(\phi(E_k),\phi(E_k))
\Big\}\\
&=&\Scal^{g}-\sum_{\ell=1}^{s}\Ric^{g}(\xi_\ell,\xi_\ell) -4ns\,,\quad\quad\quad\quad\quad\quad\quad\quad\quad\quad\quad\quad\quad\quad\quad\quad\quad\quad\quad\quad\quad\quad (\ast\ast)
\end{eqnarray*}
since  $g(E_k,E_k)=g(\phi(E_k),\phi(E_k))=1$, and the Riemannian scalar curvature is given by
\[
\Scal^{g}=\sum_{k=1}^{n}\Big\{\Ric^{g}(E_k,E_k)+\Ric^{g}(\phi(E_k),\phi(E_k))\Big\}+\sum_{\ell=1}^{s}\Ric^{g}(\xi_\ell,\xi_\ell)\,.
\]
 Now, by (\ref{not_diag}), it follows that  $\Ric^{g}(\xi_\ell, \xi_\ell)=2n\sum_{j=1}^{s}\eta_{j}(\xi_\ell)=2n$  for all $\ell\in\{1, \ldots, s\}$,  hence $\sum_{\ell=1}^{s}\Ric^{g}(\xi_\ell,\xi_\ell)=2ns$, and   the statement follows by $(\ast\ast)$.
  \pro
  \bc\label{normT2}
  Let $(M^{2n+s}, \phi, \xi_i, \eta_j, g)$  be an $\mc{S}$-manifold and let $\{e_i\}_{i=1}^{2n+s}$ be a $g$-orthonormal frame.
  The  squared norm $\|T\|^2:=\frac{1}{6}\sum_{i, j=1}^{2n+s}\|T(e_i, e_j)\|^2$ of the torsion 3-form $T$ on  $M^{2n+s}$ is given by $\|T\|^2=4ns$.
  \ec
  \pr
 This follows from  Corollary \ref{scalar_curvature} and  the  identity (see for example \cite[p.~728]{AF}).
 \[
 \Scal^{\nabla}=\Scal^{g}-\frac{3}{2}\|T\|^2\,.
 \]  
  \pro
 
\bex\label{U2exam}
Consider the Lie group $\U(2)$ with the $\mc{S}$-structure $(\phi, \xi_i, \eta_j, g)$ of CR-codimension  2 presented in \cite{TK07}.  It will be useful to adopt the notation  introduced in  \cite[Section 5.5.1]{BC} and express the distributions $\mc{D}, \mc{D}^{\perp}$ as  $\mc{D}={\rm span}\{e_1,e_2:=\phi(e_1)\}$ and $\mc{D}^{\perp}={\rm span}\{e_3=\xi_1, e_4=\xi_2\}$, respectively. 
The characteristic connection $\nabla$ on $(\U(2), \phi, \xi_i, \eta_j, g)$ was computed in \cite[Proposition.~5.2]{BC}.   Moreover,  the Ricci tensor $\Ric^{\nabla}$ of   $\nabla$  and the Riemannian Ricci tensor $\Ric^g$ 
were  described  in \cite[Proposition.~5.8]{BC}; 
\[
\Ric^\nabla=\begin{pmatrix}
-4 & 0 & 0 & 0\\
0 & -4 & 0 & 0\\
0 & 0 & 0 & 0\\
0 & 0 & 0 & 0
\end{pmatrix},\quad 
\Ric^g=\begin{pmatrix}
0 & 0 & 0 & 0\\
0 & 0 & 0 & 0\\
0 & 0 & 2 & 2\\
0 & 0 & 2 & 2
\end{pmatrix}\,.
\]
Hence, for example, for the horizontal parts  we have
$\Ric^{\nabla}(e_1, e_1)=-4$ and $\Ric^{g}(e_1, e_1)=0$, which confirm the relation $\Ric^{\nabla}(e_1, e_1)=\Ric^{g}(e_1, e_1)-2 \cdot 2 \cdot g(e_1, e_1)$, and similar for $\Ric^{\nabla}(e_2, e_2)$. Moreover, we  have $\Scal^{\nabla}=-8$ and $\Scal^{g}=4$,  and these values confirm the relation
$\Scal^{\nabla}=\Scal^{g}-6\,n\,s$, since $n=1$ and $s=2$.  Finally, Corollary \ref{normT2} yields   $\|T\|^2=8$, which confirms  part (3) of \cite[Proposition.~5.2]{BC},  where this result was obtained directly from the definition of $\|T\|^2$.
\eex

\bex\label{H3T3exam}
Consider the product  $M^{6}=H_3\times\Tg^3$, where $H_3$ is  the 3-dimensional Heisenberg group  endowed with its standard Sasakian structure
and $\Tg^3$ is the 3-dimensional torus. In Section 5.2.1 of \cite{BC} it was described an $\mc{S}$-structure $(\phi, \xi_i, \eta_j, g)$ on $M^6$ of CR-codimension 4, i.e., $s=4$ (see \cite{DL05} for the general construction).  Thus we have $n=1$,  $\mc{D}={\rm span}\{e_1, e_2\}$ and $\mc{D}^{\perp}={\rm span}\{e_3=\xi_1, \ldots, e_6=\xi_4\}$, where again we adopt the notation from  \cite{BC}.  The characteristic connection $\nabla$ on $(M^6, \phi, \xi_i, \eta_j, g)$ was described in \cite[Proposition.~5.10]{BC} and by \cite[Proposition.~5.16]{BC} we know that 
\[
\Ric^{\nabla}=\begin{pmatrix}
16 & 0 & 0 & 0 & 0 & 0\\[4pt]
0 & 16 & 0 & 0 & 0 & 0\\[4pt]
0 & 0 & 0 & 0 & 0 & 0\\[4pt]
0 & 0 & 0 & 0 & 0 & 0\\[4pt]
0 & 0 & 0 & 0 & 0 & 0\\[4pt]
0 & 0 & 0 & 0 & 0 & 0
\end{pmatrix},\quad 
\Ric^{g}=\begin{pmatrix}
24 & 0 & 0 & 0 & 0 & 0\\[4pt]
0 & 24 & 0 & 0 & 0 & 0\\[4pt]
0 & 0 & 2 & 2 & 2 & 2\\[4pt]
0 & 0 & 2 & 2 & 2 & 2\\[4pt]
0 & 0 & 2 & 2 & 2 & 2\\[4pt]
0 & 0 & 2 & 2 & 2 & 2
\end{pmatrix}\,.
\]
Hence, for example, we can confirm the relation $\Ric^{\nabla}(e_1, e_1)=\Ric^{g}(e_1, e_1)-2 \cdot 4 \cdot g(e_1, e_1)$, which takes the form  $16=24-2\cdot 4$. Similarly are confirmed the rest horizontal components. Moreover, we have $\Scal^{\nabla}=32$ and $\Scal^g=56$,
and these values confirm the relation $\Scal^{\nabla}=\Scal^{g}-6\,n\,s$, since $n=1$ and $s=4$. Corollary \ref{normT2} also yields $\|T\|^2=16$, which was obtained in \cite[Proposition.~5.10]{BC} by a different way.
\eex
 
  \bex \textnormal{(\textsf{Sasakian examples})}\label{exampleSak}
  Let $(M^{2n+1}, \phi, \xi, \eta, g)$ be a Sasakian manifold  endowed with its characteristic connection $\nabla=\nabla^{g}+\eta\wedge F$,  where $F=\frac{1}{2}\dd\eta$ is the fundamental 2-form. According to the results above we should have  
  \[
  \Ric^{\nabla}|_{\mc{D}\times\mc{D}}=\Ric^{g}|_{\mc{D}\times\mc{D}}-2\,g|_{\mc{D}\times\mc{D}}\,,\quad \Ric^{\nabla}|_{\mc{D}^{\perp}\times \mc{D}^{\perp}}=0\,,\quad \Scal^{\nabla}=\Scal^g-6\,n\,,
  \]
where $\mc{D}=\Ker(\eta)$ and $\mc{D}^{\perp}=\langle\xi\rangle$, respectively.  Below we verify these relations using previously known examples.  \\
(a) Consider  the 5-dimensional Heisenberg group $H_5$   endowed with its standard Sasakian structure $(\phi, \xi, \eta, g)$, see \cite[Section~5.2]{Pu12}, 
By the same work, it is known  that 
\[
 \Ric^{g}={\rm diag}(-2, -2, -2, -2, 4)\,,\quad \Ric^{\nabla}={\rm diag}(-4, -4, -4, -4, 0)\,.
\]  
Thus for any $i\in\{1, \ldots, 4\}$ we get $\Ric^{\nabla}(e_i, e_i)=-4$ and $\Ric^{g}(e_i, e_i)=-2$, which confirm the relation $\Ric^{\nabla}(e_i, e_i)=\Ric^{g}(e_i, e_i)-2 \cdot g(e_i, e_i)$. The scalar curvatures are  $\Scal^{\nabla}=-16$, $\Scal^{g}=-4$, hence verifying the relation $\Scal^{\nabla}=\Scal^g-6n$, as we have $n=2$. For a 5-dimensional Sasakian manifold it is also known that $\|T\|^2=8$,  which we now obtain by Corollary \ref{normT2}.\\  
(b) This is taken from \cite[Example 7.4]{FrIv}.  Consider a 4-dimensional simply connected K\"ahler-Einstein manifold $(N^{4}, g_{N}, J)$ with scalar curvature $\Scal^{g_{N}}=32$. Then there exists a circle bundle $\pi : M^{5}\to N^{4}$  over $N$ as well as a Sasakian structure $(\phi, \xi, \eta, g)$ on the total space $M^5$, such that
\[
\Ric^{g}={\rm diag}(6, 6, 6, 6,  4)\,,\quad \Ric^{\nabla}={\rm diag}(4, 4, 4, 4, 0)\,.
\]
It is now easy to see that these expressions verify  the given relations in Proposition \ref{ricp1}. 
Note also that $\Scal^{\nabla}=16=\Scal^g-6n=28-6\cdot 2$. \\
(c) Consider a 5-dimensional Sasakian manifold  $(M^{5}, \phi, \xi, \eta, g)$ endowed with its characteristic connection $\nabla=\nabla^{g}+\eta\wedge F$, such that $\Hol(\nabla)\subset\SU(2)$. Then $M$ is $\eta$-Einstein, see \cite[Proposition 4.1]{Schoe}. In particular, we have  the relations 
\[
\Ric^g=6g-2\eta\otimes\eta\,,\quad \Ric^{\nabla}=4g-4\eta\otimes\eta\,.
\]
Hence $\Ric^{g}=\Ric^{\nabla}+2g+2\eta\otimes\eta$, which yields that $\Ric^{\nabla}|_{\mc{D}\times\mc{D}}=\Ric^{g}|_{\mc{D}\times\mc{D}}-2g|_{\mc{D}\times\mc{D}}$. Obviously, we  also have $\Ric^{\nabla}|_{\mc{D}^{\perp}\times\mc{D}^{\perp}}=0$.
  \eex

 \br
On an $\mc{S}$-manifold $(M^{2n+s}, \phi,\xi_i, \eta_j, g)$    the Ricci tensor  $\Ric^{\nabla}$ associated to the characteristic connection $\nabla$
can be expressed as (see \cite[Proposition 3.1]{IP01} for the general case)   
 \[
 \Ric^{\nabla}=\Ric^{g}-\frac{1}{4}S\,, \quad S(U,V):=\sum_{a=1}^{2n+s} g\big(T(e_a,U),\,T(e_a,V)\big)\,.
 \]
 Since $\Ric^{\nabla}|_{\mc{D}\times\mc{D}}=\Ric^{g}|_{\mc{D}\times\mc{D}}-2sg|_{\mc{D}\times\mc{D}}$ , it follows that the horizontal part of  the symmetric tensor $S$ is a multiple of the Riemannian tensor $g|_{\mc{D}\times\mc{D}}$, that is,  $S|_{\mc{D}\times\mc{D}}=8sg|_{\mc{D}\times\mc{D}}$.
 According to \cite{BC}, for the examples presented above, the matrix presentation of $S$ is given
 by
 \[
 S=\begin{pmatrix}
16 & 0 & 0 & 0\\
0 & 16 & 0 & 0\\
0 & 0 & 8 & 8\\
0 & 0 & 8 & 8
\end{pmatrix}\,,\quad S=\begin{pmatrix}
32 & 0 & 0 & 0 & 0 & 0\\[4pt]
0 & 32 & 0 & 0 & 0 & 0\\[4pt]
0 & 0 & 8 & 8 & 8 & 8\\[4pt]
0 & 0 & 8 & 8 & 8 & 8\\[4pt]
0 & 0 & 8 & 8 & 8 & 8\\[4pt]
0 & 0 & 8 & 8 & 8 & 8
\end{pmatrix}
 \]
 respectively, and both these expressions confirm the relation $S|_{\mc{D}\times\mc{D}}=8sg|_{\mc{D}\times\mc{D}}$.  Similarly, for the  Sasakian examples. For instance, according to  \cite{FrIv} in case (b) of Example \ref{exampleSak}  the tensor $S$ is given by $S={\rm diag}(8, 8, 8, 8, 16)$.
 \er
 
 \subsection{Ricci-flat metric connections with parallel skew-torsion}
   Let $(M^{n}, g, T)$ be a  geometry  with parallel skew-torsion and let $\nabla=\nabla^g+\frac{1}{2}T$  be the characteristic connection.
   We say that  $(M^{n}, g, T)$ is  \textsf{$\nabla$-Einstein}  
if  $\Ric^{\nabla}$ is proportional to the metric,  
 \begin{equation}\label{nabla-einstein}
  \Ric^{g}-\frac{1}{4}S= c\cdot g\,,\quad\text{where}\quad c=\frac{\Scal^{\nabla}}{n}\,.
  \end{equation}
If a geometry with parallel skew-torsion $(M^{n}, g, T)$  is $\Ric^{\nabla}$-flat, i.e., $\Ric^{\nabla}\equiv 0$, then  it is trivially $\nabla$-Einstein with $c=\Scal^{\nabla}=0$.  If  (\ref{nabla-einstein}) is satisfied  with  $\Scal^{\nabla}\neq 0$ then  we typically  speak for \textsf{strict $\nabla$-Einstein}
 manifolds with parallel skew-torsion, see also \cite{AF14, DGP18, CGW19} for more details on $\nabla$-Einstein manifolds with (parallel) skew-torsion.    Nearly K\"ahler manifolds and nearly parallel $\G_2$-manifolds are known examples of such geometries, see for example \cite{FrIv, C16}.
\bt
An $\mc{S}$-manifold $(M^{2n+s}, \phi, \xi_i, \eta_j, g)$   is never  strict $\nabla$-Einstein.
\et
\pr 
 This follows by  Proposition \ref{ricp1}, since $\Ric^{\nabla}|_{\mc{D}^{\perp}\times \mc{D}^{\perp}}=0$.
 \pro
 However, we can prove the following interesting result which relates the $\Ric^{\nabla}$-flatness with the transverse K\"ahler-Einstein condition, 
 see  Definition \ref{transvKE}. 
 \bt\label{theorem_flat_ric}
 Let $(M^{2n+s}, \phi, \xi_i, \eta_j, g)$ be an $\mc{S}$-manifold endowed with its characteristic connection $\nabla$.
 Then $M^{2n+s}$ is $\Ric^{\nabla}$-flat if and only if $M^{2n+s}$ is transverse K\"ahler-Einstein with Einstein constant $\lambda=4s$.
 \et
\pr
According to Corollary \ref{transvRic} the Riemmanian Ricci tensor and the transverse Ricci tensor are related by
\[
\Ric^{g}(X, Y)=\Ric^{\check{g}}(X, Y)-2s\,\check{g}(X, Y)\,,
\]
for any $X, Y\in\mc{D}$, where $\check{g}$ is the transverse metric, see (\ref{TKahlerg}).
On the other hand  by Proposition \ref{ricp1}
we have  $\Ric^{\nabla}(X, Y)=\Ric^{g}(X, Y)-2\,s\,g(X, Y)$ 
for all $X, Y\in\mc{D}$. Since by (\ref{TKahlerg}) we have $\check{g}(X, Y)=g(X, Y)$ for all $X, Y\in\mc{D}$, a combination gives the relation
\begin{equation}\label{chartransverse}
\Ric^{\nabla}(X, Y)=\Ric^{\check{g}}(X, Y)-4s\,\check{g}(X, Y)
\end{equation}
for any $X, Y\in\mc{D}$. Hence $\Ric^{\nabla}|_{\mc{D}\times\mc{D}}\equiv 0$ if and only if 
$\Ric^{\check{g}}(X, Y)=4s\,\check{g}(X, Y)$ for any $X, Y\in\mc{D}$. Our claim now follows by  Proposition \ref{ricp1} since we also have 
\[
\Ric^{\nabla}(\mc{D}, \mc{D}^{\perp})=\Ric^{\nabla}(\mc{D}^{\perp}, \mc{D}^{\perp})=0\,.
\]
\pro
\bc
A Sasakian manifold $(M^{2n+1}, \phi, \xi, \eta, g)$ endowed with its characteristic connection $\nabla=\nabla^g+\eta\wedge F$ is $\Ric^{\nabla}$-flat if and only if its transverse K\"ahler structure is Einstein with Einstein constant $\lambda=4$.
\ec
\bex
Consider the  3-sphere $\Ss^3$, endowed with its standard Sasakian structure $(\phi, \xi, \eta, g)$, where  $g$ is the standard (round) metric, see also \cite[Section 2.2]{DGP18} where $g$ is denoted by $g_{-1}$. We may consider a local orthonormal frame $\{e_1, e_2, \xi\}$   such that $\phi(e_1)=e_2$, $\phi(e_2)=-e_1$, $\phi(\xi)=0$, $\eta(\xi)=1$, and $\eta(e_1)=0=\eta(e_2)$, with $\xi$ representing the Reeb vector field and $\phi$ obtained as the restriction from the (left) multiplication by   the imaginary unit $i=\sqrt{-1}$ on $\C^2$. The existence of two Ricci-flat  invariant metric connections  $\nabla^{\pm 1}$ with skew-torsion  adapted to the Sasakian structure on $\Ss^3$ is mentioned in \cite[Remark 5]{DGP18} (since  $\Ss^3\cong\SU(2)$ is isometric to a Lie group, a characteristic connection is not  unique).
The Ricci flatness of $\Ss^3$ endowed with its standard  metric and a metric connection   with parallel skew-torsion was first described within the  spinorial framework,   see  \cite[Example 7.2]{AF} and \cite[Section 5.2]{C16}.  It turns out that such a connection should  be one of the $\pm 1$-Cartan-Schouten connections, see the paragraph after Lemma 5.4 in \cite{C16}.

 Let us   explain how the Ricci flatness of $\Ss^3$ with respect to $\nabla^{\pm 1}$ follows also from Theorem \ref{theorem_flat_ric}.  It is well-known that the standard metric on $\Ss^3$ is Einstein with Einstein constant $2$, see for example \cite{C16}, hence $\eta$-Einstein with $\al=2$ and $\beta=0$ (in terms of the Definition \ref{etaE}). Then, by Theorem \ref{Main_Corres} it follows that is transverse K\"ahler Einstein with Einstein constant $\lambda=(\al+2\cdot s)=2+2\cdot 1=4$, which agrees with the desired value $4$.  
 \eex
 
Note that the construction of $\Ric^{\nabla}$-flat $\mc{S}$-geometries of higher CR-codimension $s\geq 2$ remains an interesting open problem.


\section{Holonomy features of the characteristic connection on $\mc{S}$-manifolds}\label{SubmS}

 \subsection{Holonomy features}\label{HolF}
In this section we  begin  by recalling the   notion of reducibility within the framework of geometries with parallel skew-torsion, and other relevant definitions.

\bd
A  geometry with parallel skew-torsion $(M^{n}, g,  T)$ is said to be \textsf{reducible}, if the holonomy representation of the characteristic $\nabla=\nabla^{g}+\frac{1}{2}T$  is reducible, and \textsf{irreducible} otherwise.
\ed

 \bd \label{holonomy_def}  Let $(M^{n}, g,  T)$ be a   geometry with parallel skew-torsion  such that   $TM=\mc{T}_{1}\oplus \mc{T}_2$ is an orthogonal decomposition of the tangent bundle
 for some non-trivial $\nabla$-parallel distributions $\mc{T}_{i}\subset TM$, where $\nabla=\nabla^g+\frac{1}{2}T$ is the characteristic connection. We say that $(M^{n}, g,  T)$ is \textsf{decomposable},   if    $T\in(\Lambda^3\mc{T}_1\oplus\Lambda^3\mc{T}_2)$, i.e.,  $T=T_1+T_2$ with $T_{i}\in\Lambda^3\mc{T}_{i}$. Otherwise, we say that $(M^{n}, g,  T)$ is \textsf{indecomposable}.
 \ed

In the case of reducible holonomy, and under the additional assumption that the torsion is a $\nabla$-parallel 3-form, a de Rham--type decomposition theorem holds only when the geometry is decomposable. In particular,
   \bt \textnormal{(see \cite[Lemma~3.2]{CMS21})}
 A geometry with parallel skew-torsion is decomposable if and only if it is locally isometric to a product of geometries with parallel skew-torsion.  
 \et
In the torsion-free case $(T=0)$, this recovers the local de Rham decomposition theorem: decomposability is equivalent to the reducibility of the holonomy representation.  Accordingly, a torsion-free geometry with reducible holonomy admits a local decomposition as a Riemannian product of irreducible factors corresponding precisely to the invariant summands of the holonomy representation. 
In contrast, this equivalence breaks down in the presence of torsion. Even when the holonomy representation is reducible, the torsion tensor may possess mixed components that couple the invariant subbundles, thereby obstructing any local product structure. Consequently, a geometry with skew-torsion may have reducible holonomy while nevertheless remaining locally indecomposable.

Another remarkable result concerning reducible holonomy of connections with skew-torsion is that of    Dileo and Lotta \cite{DL18}. 
 \bt \textnormal{(\cite{DL18})}
 A geometry with skew-torsion $(M^{n}, g,  T)$ (not necessarily satisfying the parallelism condition $\nabla T=0$), whose holonomy representation is reducible is   locally a  Riemannian product, provided
 that the Riemannian sectional curvature $K^g$ is non-positive. 
 \et
 
 Let us now focus on $\mc{S}$-manifolds and first point out  that there are  no $\mc{S}$-manifolds $M^{2n+s}$  of CR-codimension $s\geq 2$  of constant sectional
curvature and, for Sasakian manifolds  $(s=1)$, the unit sphere provides the unique example (see \cite[p.~210]{BG08}).
 This is due to the  fact that on an $\mc{S}$-manifold we have
 \[
 K^{g}(X, \xi_i)=1\,,\quad K^{g}(\xi_i, \xi_j)=0
 \]
 for any unit vector field $X\in\mc{D}$ and $i, j\in\{1, \ldots, s\}$ (see \cite{Blair70, CFF90}). 
 As a consequence,   $\mc{S}$-manifolds $(M^{2n+s}, \phi, \xi_i, \eta_j, g)$ cannot  satisfy the Dileo--Lotta hypothesis of everywhere non-positive sectional curvature. In particular, we see that $\mc{S}$-manifolds admit the following holonomy interpretation.  
 
\bp\label{HolS1}
An $\mc{S}$-manifold $(M^{2n+s}, \phi, \xi_i, \eta_j, g)$ endowed with its  characteristic connection $\nabla$  has reducible holonomy contained in $\U(n)$, which is indecomposable.
\ep
 \pr
 Consider the  orthogonal  decomposition 
 \[
 TM=\mc{T}_1\oplus\mc{T}_2\,,\quad \mc{T}_{1}:=\mc{D}=\Imm(\phi)=\bigcap_{i=1}^{s}\Ker(\eta_i)\,,\quad \mc{T}_2:=\mc{T}_1^{\perp}=\mc{D}^{\perp}=\Ker(\phi)=\langle\xi_1, \ldots, \xi_s\rangle\,.
 \]
The characteristic connection has parallel skew-torsion given by $T=2\bar{\eta}\wedge F$, where $\bar{\eta}=\sum_{i=1}^{s}\eta_i$.
Since $\nabla\xi_i=0$  and $\nabla\eta_i=0$ for all $i\in\{1, \ldots, s\}$, it follows that  $\mc{D}$ and $\mc{D}^{\perp}$ are $\nabla$-parallel distributions, 
hence $\nabla$ preserves both $\mc{D}$ and $\mc{D}^{\perp}$. Therefore, the holonomy representation is reducible. Moreover, the holonomy group $\Hol(\nabla)\subset\SO(2n+s)$ is contained in  the structure group  $\U(n)\times\Id_{s}$,   and  by Proposition \ref{stabilizerS} we actually know that $\fr{hol}(\nabla)\subset\fr{stab}(T)\cong\fr{u}(n)$. Recall now that the 1-forms $\eta_{j}$ are the dual 1-forms of the Reeb vector fields $\xi_i$, and the fundamental 2-form $F$ satisfies $\xi_i\lrcorner F=0$ for all $i\in\{1, \ldots, s\}$.  Hence $F$ is  a horizontal 2-form in the sense that  $F\in\Lambda^2\mc{D}$, and thus,
  the torsion 3-form $T$ is a smooth section  of the subbundle $\mc{T}_{2}\wedge\Lambda^2\mc{T}_1=\mc{D}^{\perp}\wedge\Lambda^2\mc{D}$,  that is,  
  \[
  T=2\bar{\eta}\wedge F\in \mc{D}^{\perp}\wedge\Lambda^2\mc{D}\,.
  \]
Therefore $(M^{2n+s}, g,  T)$ is indecomposable (in particular,  by Proposition \ref{2-valuedT} we have $T(X, \xi_i)=2\phi(X)\neq 0$ for any $X\in\mc{D}$ and $i\in\{1, \ldots, s\}$).
   \pro
 
 Therefore, 
 $\mc{S}$-manifolds $(M^{2n+s}, \phi, \xi_i, \eta_j, g)$ endowed with their  characteristic connection $\nabla=\nabla^g+\bar{\eta}\wedge F$  provide examples of  geometries with parallel skew-torsion with  reducible but indecomposable holonomy, for  any CR-dimension $n\geq 1$ and CR-codimension $s\geq 1$. 
Of course, this includes the better understood odd-dimensional Sasakian case $(s=1)$, 
where  the torsion has trivial kernel, $\Ker(T)=\{0\}$.  However, the class of $\mc{S}$-manifolds is considerably broader than the Sasakian one, since it contains both odd- and even-dimensional manifolds.  Such a representative example is the Lie group $\U(2)$  mentioned in the previous section.

\br\label{RiemannS} 
Proposition \ref{HolS1} illustrates the different behaviour of   reducible holonomy in the presence of torsion, compared to the torsion-free scenario ($T=0$).
  Indeed, according to \cite[Theorem 2.1]{TP13}   any $\mc{S}$-manifold $(M^{2n+s}, \phi, \xi_i, \eta_j, g)$ is locally a Riemannian product of a  $\sqrt{s}$-Sasakian manifold and a  flat manifold of dimension  $s-1$.
We mention that this $(s-1)$-dimensional flat factor  corresponds exactly to (the fibers of) $\Ker(T)$.
Thus, the torsion tensor of the characteristic connection can be captured in the splitting phenomenon that appears in the torsion-free setting, while at the same time couples the horizontal and vertical distributions so that we get an indecomposable geometry with parallel skew-torsion.
\er


\subsection{Locally defined submersions  induced by the characteristic connection}
Let $(M^{n}, g, T)$ be a geometry with parallel skew-torsion. 
Recall that the holonomy algebra $\fr{hol}(\nabla)\subset\fr{so}(n)$ of the characteristic connection $\nabla=\nabla^g+\frac{1}{2}T$  acts naturally on $TM$. 
 In the sequel we will assume that  $TM$ decomposes orthogonally  into $\nabla$-parallel distributions, as in Section \ref{HolF},  which here it is more convenient to  denote by $\mc{H}, \mc{V}$ respectively,  i.e.,  $TM=\mc{H}\oplus\mc{V}$.  
 \bd
The splitting $TM=\mc{H}\oplus\mc{V}$ is said to be \textsf{admissible} if the projection of the torsion 3-form $T$ of the characteristic connection $\nabla$  onto the summand $\mc{H}\otimes\Lambda^2\mc{V}$ vanishes,
\[
0={\rm pr}_{\mc{H}\otimes\Lambda^2\mc{V}}T\in \mc{H}\otimes\Lambda^2\mc{V}\,. 
\]
In this case $\mc{H}$ (respectively $\mc{V}$) is called the horizontal (respectively, vertical) distribution of the admissible 
splitting.
\ed
Let $(M^{n}, g, T)$ be a geometry with parallel skew-torsion and admissible splitting $TM=\mc{H}\oplus\mc{V}$, as defined above. Since  the torsion $T$ should vanish on $\mc{H}\otimes\Lambda^2\mc{V}$, it decomposes into the remaining three parts of the decomposition of $\Lambda^3(\mc{H}\oplus\mc{V})$, 
\[
T=T^{\mc{H}}+T^{\rm mix}+T^{\mc{V}}\,,
\]
with $T^{\mc{H}}:={\rm pr}_{\Lambda^3\mc{H}}T\in\Lambda^3\mc{H}$,  $T^{\rm mix}:={\rm pr}_{\mc{V}\otimes \Lambda^2\mc{H}}T\in \mc{V}\otimes \Lambda^2\mc{H}$ and  $T^{\mc{V}}:={\rm pr}_{\Lambda^3\mc{V}}T \in\Lambda^3\mc{V}$. 
Obviously, if a   geometry  with parallel skew-torsion $(M^{n}, g, T)$   is decomposable  in terms of the Definition \ref{holonomy_def}, then $T^{\rm mix}=0$.  Moreover, it is easy to see that $\mc{V}$ is an integrable distribution and by \cite[Remark 3.15]{CMS21} it is known that with any admissible splitting of $TM$ one can associate a locally defined submersion $\pi : M\to N$ enjoying many important  features. We list  some of them in the form of a theorem. Here we mainly follow  \cite[Section 3]{MS24}. 
 \bt\label{stCMS}\textnormal{(\cite{CMS21, MS24})}
 Suppose that  $(M^{n}, g, T)$ is a geometry with parallel skew-torsion and an admissible splitting $TM=\mc{H}\oplus\mc{V}$, as defined above. This  induces a {\it locally defined} Riemannian submersion $\pi : (M, g)\to (N, g')$,  having the following properties:
\begin{itemize}
\item[(1)] $N$ is the leaf space of the integrable distribution $\mc{V}$ and in particular  $\mc{V}$ is   the vertical distribution of $\pi$. Moreover,    $\pi$ has totally geodesic fibers.
\item[(2)] The horizontal part $T^{\mc{H}}$ is projectable, i.e., there exists a 3-form $T'\in\Lambda^3TN$ on $N$ such that $T^{\mc{H}}=\pi^* T'$. 
\item[(3)] $(N, g', T')$ is  a geometry with parallel skew-torsion  with characteristic connection given by $\nabla'=\nabla^{g'}+\frac{1}{2}T'$. 
\item[(4)] The characteristic connections $\nabla$ and  $\nabla'$ are $\pi$-related, 
\[
\nabla'_{X}Y=\pi_{*}(\nabla_{X^{h}}Y^{h})\,,\quad \forall \ X, Y\in TN\,, 
\]
where $X^{h}, Y^{h}\in TM$ are the horizontal lifts of $ X, Y\in TN$ {\rm(}the Riemmanian connections $\nabla^g$ and $\nabla^{g'}$ are also $\pi$-related\,{\rm)}.
\end{itemize}
 \et
 \br
 (a) We emphasize  that this result is only of  local 
nature: one restricts to a sufficiently small  open set in $M$,  such that $N$, i.e., the leaf space of the  foliation $\scF_{\mc{V}}$ corresponding to  $\mc{V}=\mc{D}^{\perp}$ is   smooth, as the quotient of this open set by $\scF_{\mc{V}}$. \\
 (b) The holonomy algebra $\fr{hol}(\nabla)$ of  the characteristic connection on a geometry with parallel skew-torsion $(M^{n}, g, T)$ can be hard to compute.
However, recall by Section \ref{MCST} that $\fr{hol}(\nabla)$  is contained in the stabilizer algebra $\fr{stab}(T)$ of the torsion 3-form $T$,  i.e., $\fr{hol}(\nabla)\subset\fr{stab}(T)$. 
 In general, these   subalgebras  do not coincide and one can consider  the notion
of   \textsf{canonical $\fr{g}$-splittings}  for some intermediate algebra $\fr{g}$, i.e., $\fr{hol}(\nabla)\subset\fr{g}\subset\fr{stab}(T)$.\footnote{In this note by the inclusion $\subset$ we generally mean $\subseteq$.}
\er 

\bd \textnormal{(\cite[Definition~3.3]{MS24})}
Let $\fr{g}\subset\fr{so}(n)$ be a Lie algebra such that $\fr{hol}(\nabla)\subset\fr{g}\subset\fr{stab}(T)$. The representation of $\fr{g}$ on $\R^n$ decomposes into an orthogonal sum of irreducible modules $\fr{h}_{a}$ and $\fr{v}_{j}$ such that $\fr{so}(\fr{h}_{a})\cap\fr{g}\neq 0$ for all $a$ and $\fr{so}(\fr{v}_j)\cap\fr{g}=0$ for all $j$. Then, the \textsf{canonical $\fr{g}$-splitting} of $TM$ is defined by taking $\mc{H}$ to be the subbundle associated to $\fr{h}:=\oplus_{a}\fr{h}_a$ and $\mc{V}$ to be the subbundle associated to $\fr{v}:=\oplus_{j}\fr{v}_{j}$.
\ed
Obviously,   in this case the torsion 3-form $T$ corresponds to a $\fr{g}$-invariant  element in $\Lambda^3\R^n$, i.e.,   $T\in (\Lambda^3\R^n)^{\fr{g}}$.
Moreover, since $\fr{hol}(\nabla)\subset\fr{g}$ the splitting $TM=\mc{H}\oplus\mc{V}$ associated to $\R^{n}=\fr{h}\oplus\fr{v}$ is $\nabla$-parallel.   In fact,  any canonical $\fr{g}$-splitting is admissible (\cite[Lemma 3.4]{MS24}), thus  it locally defines a Riemannian submersion $\pi : (M, g)\to (N, g')$, enjoying the properties listed in Theorem \ref{stCMS}.   This is  called the \textsf{canonical $\fr{g}$-submersion}.     The canonical $\fr{g}$-submersion with  $\fr{g}$ being the holonomy algebra, i.e., $\fr{g}=\fr{hol}(\nabla)$, yields the notion of the  standard submersion introduced in \cite{CMS21}. 
 In \cite{CMS21} the authors also introduced the notion of geometries with parallel skew-torsion of \textsf{special type}, see also \cite[Section 7]{MS24}.
 \bd
A geometry with parallel skew-torsion $(M^{n}, g, T)$ is called of \textsf{special type} when the holonomy algebra  $\fr{hol}(\nabla)$ of the characteristic connection $\nabla$ acts trivially on the vertical distribution $\mc{V}\neq 0$ of the standard submersion.
 \ed
Observe that $(M^{n}, g, T)$ is of special type  if and only if $\mc{V}$ is spanned by $\nabla$-parallel vector fields. Sasakian manifolds endowed with their characteristic connection (\cite{FrIv}) form one of the most classical classes of geometries with parallel skew-torsion of special type. In this case  $\mc{V}$ is spanned by the Reeb vector field $\xi$. Obviously, their higher-dimensional analogues, $\mc{S}$-manifolds, likewise furnish many   new examples of geometries with parallel skew-torsion of special type.


 \subsection{The canonical $\fr{u}(n)$-submersion of  $\mc{S}$-manifolds}
 Let  $(M^{2n+s}, \phi, \xi_i, \eta_j, g)$ be an $\mc{S}$-manifold of CR-dimension $n\geq 1$ and CR-codimension $s\geq 1$.  By Proposition \ref{stabilizerS}  for the torsion 3-form $T$ of the characteristic connection $\nabla$ we have $\fr{stab}(T)\cong\fr{u}(n)=\Lie(\U(n))$ and hence $\fr{hol}(\nabla)\subset\fr{u}(n)$. Since $\fr{hol}(\nabla)$ can be in general unknown, next we will consider the canonical $\fr{g}$-splitting  induced by  $\fr{g}=\fr{stab}(T)\cong\fr{u}(n)$.  An application   of  Theorem \ref{stCMS} gives that

\bt\label{standSm}
Let $(M^{2n+s}, \phi, \xi_i, \eta_j, g)$ be an $\mc{S}$-manifold endowed with the characteristic connection $\nabla$ with parallel skew-torsion $T=\sum_{i=1}^{s}\eta_i\wedge\dd\eta_i=2\bar{\eta}\wedge F$. Then the following hold:\\
{\rm (1)} $TM$ admits a canonical $\fr{u}(n)$-splitting given by  
\begin{equation}\label{admis}
TM=\mc{H}\oplus\mc{V}\,,\quad \mc{H}=\mc{D}=\Imm(\phi)\,,\quad \mc{V}=\mc{D}^{\perp}=\Ker(\phi)\,.
\end{equation}
The components of the torsion 3-form $T$  with respect to this admissible splitting are such that 
\[
T^{\mc{H}}= 0=T^{\mc{V}}\,, \quad T^{\rm mix}={\rm pr}_{\mc{V}\otimes\Lambda^2\mc{H}}T=\sum_{i=1}^{s}\eta_i\wedge\dd\eta_i=T\,.
\]
Thus, there is locally defined  Riemannian submersion $\pi : M^{2n+s}\to N^{2n}$, where $N^{2n}$ is the leaf space  of the foliation $\scF_{\mc{V}}$ generated by the integrable vertical distribution $\mc{V}$ (and hence of the  characteristic foliation of  $M^{2n+s}$).\\ 
{\rm(2)} The 3-form $T'\in\Lambda^3TN$ on the base space satisfying the relation $T^{\mc{H}}=\pi^*T'$  vanishes, and thus   the characteristic connections $\nabla, \nabla'$ on $M, N$, respectively,  satisfy 
the relation
\[
\nabla^{g'}_{X}Y=\nabla'_{X}Y=\pi_{*}(\nabla_{X^{h}}Y^{h})\,,\quad X, Y\in TN\,,
\]
where $\nabla^{g'}$ is the Levi-Civita connection of the base $(N, g')$ and  $X^{h}$ is the horizontal lift of $X\in TN$. \\
{\rm(3)} The base space $N$ admits a K\"ahler structure $(g', \mathbb{J})$ where the   complex structure $\mathbb{J} $ is defined by
\begin{equation}\label{acJa}
\mathbb{J}X=\pi_{*}\phi(X^{h})\,,\quad X\in TN\,.
\end{equation}
\noindent {\rm(4)}  If $(M^{2n+s}, \phi, \xi_i, \eta_j, g)$  is a compact  regular $\mc{S}$-manifold, then  the  canonical $\fr{u}(n)$-submersion $\pi : M^{2n+s}\to N^{2n}$ is globally defined and  coincides with the $\BLY$-submersion, i.e., $\pi=\pi_s$ for all $s\geq 1$.
\et
\pr
 (1)
We already know that   both the distributions 
\[
\mc{D}=\Imm(\phi)=\bigcap_{i=1}^{s}\Ker(\eta_i)\,,\quad \mc{D}^{\perp}=\Ker(\phi)=\langle\xi_1, \ldots, \xi_s\rangle
\]
 are $\nabla$-parallel distributions, hence they are both $\fr{hol}(\nabla)$-invariant.
The splitting $TM=\mc{H}\oplus\mc{V}=\Imm(\phi)\oplus\Ker(\phi)$  described in (\ref{admis})
 is the canonical $\fr{u}(n)$-splitting, since the horizontal and vertical distributions   are   $\fr{stab}(T)$-invariant, as well. This is because  $\fr{u}(n)$ acts in the standard way on $\mc{D}$ via the almost complex structure $\phi|_{\mc{D}}$ and trivially on $\mc{D}^{\perp}$.
Actually, the torsion 3-form $T$ is  $\fr{u}(n)$-invariant.  Indeed,  with respect to the action of $\fr{stab}(T)\cong\fr{u}(n)$ it is easy to see that $F\in(\Lambda^{2}\mc{H})^{\fr{u}(n)}$, hence we obtain the relation $T\in\mc{V}\otimes(\Lambda^2\mc{H})^{\fr{u}(n)}$. 
Thus we get an admissible splitting  with integrable vertical distribution $\mc{V}$ generated by the characteristic vector fields $\xi_1, \ldots, \xi_s$, and horizontal distribution $\mc{H}=\mc{D}=\Imm(\phi)$. 
  Theorem \ref{stCMS} then guarantees the existence of a  locally defined Riemannian submersion  $\pi : M^{2n+s}\to N^{2n}$ associated to  this admissible splitting, which is the canonical $\fr{g}$-submersion  for $\fr{g}=\fr{u}(n)$. We finally observe that  the foliation $\scF_{\mc{V}}$  generated by the vertical distribution $\mc{V}$ must coincide with the characteristic foliation $\scF$, since $\mc{V}=\mc{D}^{\perp}$. This proves (1).\\
 (2) Since $T^{\mc{H}}$ vanishes identically on $M^{2n+s}$, we have $\pi^*T'=0$ in terms of part (2) in  Theorem \ref{stCMS}, which implies that $T'=0$.  Thus, according to part (3) of Theorem \ref{stCMS} we get $\nabla'=\nabla^{g'}$  and the result follows by  part (4) of Theorem \ref{stCMS}.\\
(3)   We have $\mc{L}_{\xi_i}\phi=0$ for all $i\in\{1, \ldots, s\}$,  hence obviously 
 $\phi$ projects to the almost complex structure $\Ja$ on the $2n$-dimensional leaf space, 
\[
\Ja^2X=\pi_{*}\big(\phi^2(X^{h})\big)=\pi_{*}(-X^{h})=-\pi_{*}(X^{h})=-X\,.
\]
For the integrability of $\Ja$, it is sufficient to  show the $\nabla^{g'}$-parallelism of $\Ja$, $\nabla^{g'}\Ja=0$.
By the previous part, we have  that $\nabla^{g'}=\nabla'$ and the characteristic connections  $\nabla$ on $M$ and $\nabla'$ on $N$  are $\pi$-related. However  $\Ja$ is essentially the restriction of $\phi$ on $\mc{H}=\mc{D}$, and  by Theorem \ref{Torsion_skewS} the endomorphsim $\phi$ is $\nabla$-parallel.
Thus, since $\nabla, \nabla'$ are $\pi$-related, $\Ja$ is also $\nabla'$-parallel, and our claim follows by the identification  $\nabla^{g'}=\nabla'$.  An alternative is based on the direct computation of the  Nijenhuis  tensor $N_{\Ja}$ of $\Ja$. This essentially follows  the method presented in \cite{BLY73},  and we avoid to present the details.  \\ 
\noindent Now,  $\mc{V}$ is generated by Killing vector fields and the $\phi$-compatible  Riemannian metric $g$  on $M^{2n+s}$ projects to a metric  $g'$ on $N$, defined by
\[
g'(X, Y)\circ\pi=g(X^{h}, Y^{h})\,,\quad \forall\ X, Y\in TN\,.
\]
This is the metric induced on $N$ by the canonical $\fr{u}(n)$-submersion $\pi : M^{2n+s}\to N^{2n}$. It is easy to see that this is a Hermitian metric, thus $(N, g', \Ja)$ is a Hermitian manifold. It is K\"ahler, since the 2-form $\Om$ defined by $\Om(X, Y)=g'(X, \Ja Y)$ satisfies $\pi^*\Om=F$ and $\dd F=0$ implies $\dd\Om=0$ since $\pi^*$ is an injection.\\
(4) Under the assumptions in the statement,   the characteristic vector fields $\xi_1, \ldots, \xi_s$ generate an action of the $s$-dimensional torus $\Tg^s$ on the total space $M^{2n+s}$ and the leaves of the foliation $\scF_{\mc{V}}$ associated to the vertical distribution $\mc{V}$  
of  the  canonical $\fr{u}(n)$-submersion $\pi$ are diffeomorphic  to $\Tg^s$. Thus, in this case  $\pi$ becomes a principal $\Tg^s$-bundle over $N=M^{2n+s}/\scF_{\mc{V}}$. The claim now  follows by the identifications  $\mc{H}=\scH^{2n}$, $\mc{V}=\scV^{s}$, where $\scH^{2n}$ and $\scV^{s}$ are described  in (\ref{vh}), and the statement in  Theorem \ref{BYL}.  \pro

We are now able to summarize as follows:  
\bc
Any $\mc{S}$-manifold admits a locally defined Riemannian submersion over a K\"ahler manifold and is a geometry with parallel skew-torsion of special type.
\ec
Note that for $s=1$ this is a well-known statement (see \cite[Section 7.1]{MS24}).  

\subsection{A characterization of $\mc{S}$-manifolds}
Let $(M, g, T)$ be a geometry with parallel skew-torsion and consider the canonical $\fr{g}$-splitting with    $\fr{g}=\fr{stab}(T)$, 
\[
TM=\mc{H}\oplus\mc{V}\,,
\]
as defined above. 
In \cite{MS24} the authors studied the  case where  the vertical space $\mc{V}$ is 1-dimensional, and the horizontal space $\mc{H}$ is an irreducible representation of $\fr{stab}(T)$. They referred to this situation as \textsf{almost irreducible case}. Here we need to  generalize this definition. 

\bd\label{stabTHM}
Let  $(M^{n}, g, T)$  be a geometry with parallel skew-torsion  admitting a canonical  $\fr{stab}(T)$-invariant splitting $TM=\mc{H}\oplus\mc{V}$. We say that $\fr{stab}(T)$ acts \textsf{$s$-irreducibly} on   $TM$ if the vertical distribution $\mc{V}$ is $s$-dimensional and $\mc{H}$ is an irreducible representation of $\fr{stab}(T)$. In this case 
  $(M^{n}, g, T)$ is said to be  a   \textsf{holonomy $s$-irreducible} geometry with parallel skew-torsion.
\ed
Obviously, holonomy $1$-irreducible  is the  almost irreducible case of \cite{MS24}, with Sasakian manifolds being the most characteristic examples.
We will prove the following generalization of \cite[Theorem~7.1]{MS24}.   
 
 \bt\label{SmnfdsCHAR}
Let $(M, g, T)$ be a geometry with parallel skew-torsion such that $\fr{stab}(T)$ acts $s$-irreducibly on $TM$, with canonical splitting $TM=\mc{H}\oplus\mc{V}$. Moreover, assume that  $(M, g, T)$ has no local $s$-dimensional factor, and that  $T^{\mc{V}}=0$.  Then $(M, g, T)$ is an $\mc{S}$-manifold.
 \et
 \pr
Let $\nabla=\nabla^g+\frac{1}{2}T$ be the characteristic connection on $(M, g, T)$.  
We may assume that the vertical distribution  $\mc{V}$ in the  canonical $\fr{stab}(T)$-splitting of $TM$  is generated by  $\nabla$-parallel vector fields 
$\xi_1, \ldots, \xi_s\in TM$, i.e., $\mc{V}=\langle\xi_1, \ldots, \xi_s\rangle$ with $\nabla\xi_i=0$ for all $i\in\{1, \ldots, s\}$.  
 The horizontal part $T^{\mc{H}}$ of $T$ must vanish, see  \cite[Lemma 3.7]{MS24}, while by assumption  $T^{\mc{V}}=0$.  Thus, the torsion 3-form $T\in(\Lambda^3\R^n)^{\fr{stab}(T)}$ must lie in 
 \[
 \mc{V}\otimes(\Lambda^{2}\mc{H})^{\fr{stab}(T)}\,.
 \]
Since by our assumption $\mc{H}$ is irreducible and  $(M, g, T)$ has no $s$-dimensional local factor, we thus conclude that  our geometry $(M, g, T)$ is  indecomposable. Therefore we may assume, without loss of generality, that
\begin{equation}\label{supT}
T=2\sum_{i=1}^{s}\eta_{i}\wedge F\,,
\end{equation}
for some  horizontal 2-form $F\in\Lambda^2\mc{H}$ which is $\fr{stab}(T)$-invariant,  where $\eta_i$ are the dual 1-forms of $\xi_i$, i.e., $\eta_i(\xi_j)=\delta_{ij}$ for all $i, j\in\{1, \ldots, s\}$.   Hence we have $\xi_i\lrcorner F=0$  for all $i\in\{1, \ldots, s\}$.  Consider now the bundle endomorphism $\phi : TM\to TM$ defined by
\[
\phi(X)=-(\iota_{X}F)^{\sharp}
\]
where $\sharp : T^*M\to TM$ is the usual isomorphism   $T^*M\ni\al\longmapsto \al^{\sharp}\in TM$ 
with $g(\al^{\sharp}, Y)=\al(Y)$. 
 Obviously,   $\phi$ is a $\fr{stab}(T)$-invariant endomorphism satisfying  $\phi(\xi_i)=0$ for all $i$, and  
\[
g(X, \phi(Y))=-(\iota_{Y}F)(X)=-F(Y, X)=F(X, Y)
\]
  for all $X, Y\in\Gamma(TM)$. 
It follows that $\phi$ is skew-adjoint with respect to $g$, and moreover that $\phi^2 : TM\to TM$ is $\fr{stab}(T)$-invariant and  self-adjoint with $\phi^2(\xi_i)=0$, for all $i$.   Therefore,  on the $\fr{stab}(T)$-irreducible horizontal distribution $\mc{H}$  by Schur's lemma we  must have $\phi^2=-\lambda\Id_{\mc{H}}$ 
 for some $\lambda>0$.  After rescaling $g$ and $\xi_i$, we may assume that $\lambda=1$ (see also \cite{MS24} for $s=1$ and similar arguments). Consequently,  $\phi^2$  is given by
  \[
 \phi^2=-\Id+\sum_{i=1}^{s}\eta_{i}\otimes\xi_i\,,
 \]
and we can identify $\mc{H}=\Imm(\phi)=\bigcap_{i=1}^{s}\Ker(\eta_i)$, such that $TM=\mc{H}\oplus\mc{V}=\Imm(\phi)\oplus\Ker(\phi)$.
Moreover, the invariant complex structure $\phi|_{\mc{H}}$ allows to view $\mc{H}$ as a complex representation of $\fr{stab}(T)$. Thus, $M$  must be $(2n+s)$-dimensional with $2n=\dim\mc{H}$, and the family  $(\phi, \xi_i, \eta_j, g)$ defines a metric $f$-structure on $M$ with fundamental 2-form $F$.\\
\noindent Let us    use the torsion 3-form $T$ to prove that $2F=\dd\eta_i$ for all $i\in\{1, \ldots, s\}$, which means that actually   $(M^{2n+s}, \phi, \xi_i, \eta_j, g)$ is an almost $\mc{S}$-manifold. 
Based on (\ref{supT}),  a direct computation (as in Proposition \ref{2-valuedT}) shows that
\[
T(X, Y)=2\sum_{i=1}^{2}\{F(X, Y)\xi_i-\eta_{i}(X)\phi(Y)+\eta_{i}(Y)\phi(X)\}
\]
for all $X, Y\in TM$. 
This gives $T(\xi_i, \xi_j)=0$, and hence $[\xi_i, \xi_j]=0$ for all $i, j\in\{1, \ldots, s\}$. Moreover,   we get $T(X, \xi_i)=2\phi(X)$ and 
   $T(X, Y, \xi_i)=T(\xi_i, X, Y)=2F(X, Y)$, for all  $X, Y\in\mc{H}$, that is, $\xi_i\lrcorner T=2F$ for all $i\in\{1, \ldots, s\}$. 
   Then we also see that
   \begin{equation}\label{cond1}
   \nabla^{g}_{X}\xi_i=-\frac{1}{2}T(X, \xi_i)=-\phi(X)\,,\quad T(X, Y)=2F(X, Y)\bar{\xi}, \quad \forall \ X, Y\in\mc{H}\,,
   \end{equation}
   where $\bar{\xi}=\sum_{i=1}^{s}\xi_i$. 
Now, because the conditions $\nabla\xi_i=0$ and $\nabla\eta_i=0$ are equivalent, it follows that $(\nabla^{g}_{X}\eta_i)Y=\frac{1}{2}T(X, Y, \xi_i)=F(X, Y)$ for all $X, Y\in\mc{H}$. Hence
\[
\dd\eta_i(X, Y)=(\nabla^g_{X}\eta_i)Y-(\nabla^g_{Y}\eta_i)Y=2F(X, Y)\,.
\]
We also get $\dd\eta_i(X, \xi_j)=0=\dd\eta_{i}(\xi_j, \xi_k)$ for all $i, j, k\in\{1, \ldots, s\}$, and we deduce that 
 $\dd\eta_i=2F$ for all $i\in\{1, \ldots, s\}$.  Thus $(M^{2n+s}, \phi, \xi_i, \eta_j, g)$  is an almost $\mc{S}$-manifold. 
 Now our conclusion follows by Theorem \ref{Torsion_skewS}; our manifold $(M, g, T)$ is a geometry with parallel skew-torsion, and by Theorem \ref{Torsion_skewS} an almost $\mc{S}$-manifold $(M^{2n+s}, \phi, \xi_i, \eta_j, g)$ 
 admits a  metric connection with skew-torsion preserving the structure if and only if $(M^{2n+s}, \phi, \xi_i, \eta_j, g)$ is an $\mc{S}$-manifold, and in this case the torsion is $\nabla$-parallel.
  \pro
 
 
 \br
 For $s=1$ the condition $T^{\mc{V}}=0$ is trivially satisfied,   
 and  Theorem \ref{SmnfdsCHAR} reduces to \cite[Thm.~7.1]{MS24}. We   also observe that    the proof in \cite[Thm.~7.1]{MS24} does not explicitly rely on the fact that for $s=1$ such a triple $(M,g,T)$ is a contact metric manifold, and hence  the conclusion follows directly from \cite[Theorem 8.4]{FrIv}. Instead, the authors proceed by showing that the manifold $M$ satisfies the Sasakian identities.

Of course, such an approach also applies to our more general  scenario.  For instance, in the proof above we concluded that $(M, g, T)$ is an almost $\mc{S}$-manifold  $(M^{2n+s}, \phi, \xi_i, \eta_j, g)$.  
   Hence, to prove that $(M^{2n+s}, \phi, \xi_i, \eta_j, g)$ is an $\mc{S}$-manifold, one can either show that  $M^{2n+s}$  is normal or, for example,  that its Riemannian curvature tensor satisfies
  \begin{equation}\label{RicgXY}
  R^{g}(X, Y)\xi_i=\sum_{j=1}^{s}\{\eta_{j}(X)\phi^2(Y)-\eta_{j}(Y)\phi^2(X)\}
  \end{equation}
  for all $X, Y\in TM$ and $i\in\{1, \ldots, s\}$, see \cite[Theorem 4.3]{CFF90}. Let us employ the second criterion. \\ 
\noindent  Choose $X, Y\in\mc{H}$ and $Z=\xi_i$ for some $i\in\{1, \ldots, s\}$. Then (\ref{RicgXY}) gives $R^{g}(X, Y)\xi_i=0$ for all $i\in\{1, \ldots, s\}$, which we need to confirm. Indeed, a long but direct computation shows that
 \begin{equation}\label{3-curvature}
 R^{\nabla}(X, Y)Z=R^{g}(X, Y)Z+\frac{1}{4}T(T(Y, Z), X)-\frac{1}{4}T(T(X, Z), Y)+\frac{1}{2}T(T(X, Y), Z)
 \end{equation}
 for all $X, Y, Z\in TM$.  Now, the condition  $\nabla\xi_i=0$ for all $i\in\{1, \ldots, s\}$   implies $R^{\nabla}(X, Y)\xi_i=0$ for all such $i$ and    $X, Y\in\mc{H}$. Thus,  by (\ref{3-curvature}) we obtain
 \begin{eqnarray*}
 R^{g}(X, Y)\xi_i&=&\frac{1}{4}\big\{T(T(X, \xi_i), Y)-T(T(Y, \xi_i), X)\big\}-\frac{1}{2}T(T(X, Y), \xi_i)\\
 &\overset{(\ref{cond1})}{=}&\frac{1}{2}\big\{T(\phi(X), Y)-T(\phi(Y), X)\big\}-\frac{2}{2}F(X, Y)T(\bar{\xi}, \xi_i)\\
 &=&\big\{F(\phi(X), Y)-F(\phi(Y), X)\big\}\bar{\xi}=0
 \end{eqnarray*}
 since $T(\xi_i, \xi_j)=0$ for all $i, j\in\{1, \ldots, s\}$ and hence $T(\bar{\xi}, \xi_i)=0$ for all $i\in\{1, \ldots, s\}$,  
 and moreover $F(\phi(X), Y)+F(X, \phi(Y))=0$ for all $X, Y\in TM$, see for example (2.2) in \cite{BC}.\\ 
\noindent Similarly, for $X\in\mc{H}$, $Y=\xi_k$ and $Z=\xi_i$ for some $k, i\in\{1, \ldots, s\}$, by (\ref{RicgXY}) we should have $R^{g}(X, \xi_k)\xi_i=-\phi^2(X)=X$.  To confirm this, Recall by Corollary \ref{curvatCorl} that  $R^{\nabla}(X, \xi_k)\xi_i=0$, thus equation (\ref{3-curvature}) gives
 \begin{eqnarray*}
R^{g}(X, \xi_k)\xi_i&=&\frac{1}{4}\big\{T(T(X, \xi_i), \xi_k)-T(T(\xi_k, \xi_i), X)\big\}-\frac{1}{2}T(T(X, \xi_k), \xi_i)\\
&=&\frac{1}{2}T(\phi(X), \xi_k)-T(\phi(X), \xi_i)=\phi^2(X)-2\phi^2(X)=-X+2X=X\,.
 \end{eqnarray*}
 Similarly is treated the case where $X, Y\in\mc{V}$.
\er


\subsection{The curvature relations of the canonical submersion}
We will now focus on the curvature relations  of the canonical $\fr{u}(n)$-submersion  $\pi : M^{2n+s}\to N^{2n}$ 
 of an $\mc{S}$-manifold $(M^{2n+s}, \phi, \xi_i, \eta_j, g)$, as described above.  From now on, for simplicity, we will refer to this locally defined Riemannian submersion by the term \textsf{canonical submersion} of $M$.   Let $\mc{H}=\mc{D}$ and $\mc{V}=\mc{D}^{\perp}$ be the corresponding horizontal and vertical distributions, as in (\ref{admis}). 
 The fundamental tensor fields of $\pi$ are the O'Neill's tensors $\mc{A}$ and $\mc{T}$, which are defined by
\begin{eqnarray*}
\mc{A}_{X}Y&=&(\nabla^{g}_{X_{\mc{H}}}Y_{\mc{H}})_{\mc{V}}+(\nabla^g_{X_{\mc{H}}}Y_{\mc{V}})_{\mc{H}}\,,\\
\mc{T}_{X}Y&=&(\nabla^{g}_{X_{\mc{V}}}Y_{\mc{H}})_{\mc{V}}+(\nabla^{g}_{X_{\mc{V}}}Y_{\mc{V}})_{\mc{H}}\,,
\end{eqnarray*}
where  $X_{\mc{H}}$ and $X_{\mc{V}}$ denote the   projections of a vector field $X\in TM$ on $\mc{H}$ and $\mc{V}$, respectively. 
We already know that  $\mc{T}$ vanishes identically, since the fibers of $\pi$ are totally geodesic (and can be directly verified in our framework). 
Regarding the  tensor  $\mc{A}$,  we have the relation $\mc{A}_{X} = -\frac{1}{2}T^{\rm mix}_{X}=-\frac{1}{2}(X\lrcorner T^{\rm mix})$ (see   \cite{MS24}),   which for our case, because of the identification $T^{\rm mix}=T$,  reduces to $\mc{A}_{X}=-\frac{1}{2}(X\lrcorner T)$. In particular, 
\bl\label{previous}
{\rm(1)} The tensor field $\mc{A}$ satisfies the relation 
\begin{eqnarray*}
g(A_{X}Y, Z)&=&-\frac{1}{2}\Big\{T(X_{\mc{H}}, Y_{\mc{H}}, Z_{\mc{V}})+T(X_{\mc{H}}, Y_{\mc{V}}, Z_{\mc{H}})\Big\}\\
&=&-F(X_{\mc{H}}, Y_{\mc{H}})\sum_{j=1}^{s}\eta_{j}(Z_{\mc{V}})+F(X_{\mc{H}}, Z_{\mc{H}})\sum_{j=1}^{s}\eta_{j}(Y_{\mc{V}})\,.
\end{eqnarray*}
Hence,   $\mc{A}_{X}Y=-F(X, Y)\bar{\xi}=-\frac{1}{2}T(X, Y)$ and    $\mc{A}_{X}\xi_k=-\phi(X)=-\frac{1}{2}T(X, \xi_k)$, 
for all $X, Y\in\mc{H}$  and $k\in\{1, \ldots, s\}$. \\
{\rm (2)} For any vertical vector field $V\in\mc{V}$ and for any basic vector field $Z\in\mc{H}$ the characteristic connection $\nabla$ on $M^{2n+s}$  is such that
 \[
 (\nabla_{V}Z)_{\mc{H}}=-2\sum_{j=1}^{s}\eta_{j}(V)\phi(Z)\,.
 \]
{\rm (3)} Let   $
 \{e_{a}\}_{a=1}^{2n+s}=\big\{E_1, \ldots, E_n, \phi(E_1), \ldots, \phi(E_n), \xi_1, \ldots, \xi_s\big\}$
  be a local $g$-orthonormal $\phi$-adapted frame of $M^{2n+s}$.  The quantity   $g(\mc{A}_X, \mc{A}_Y):=\sum_{a=1}^{2n}g(\mc{A}_{X}e_{a}, \mc{A}_{Y}e_{a})$ 
    is a multiple of $g(X, Y)$ for any $X, Y\in\mc{H}$, in particular $ g(\mc{A}_X, \mc{A}_Y)=s\,g(X, Y)$. 
Thus, we also have
\[
\sum_{a, b=1}^{2n}g(\mc{A}_{e_{a}}e_{b}, \mc{A}_{e_{a}}e_{b})=s\sum_{a, b=1}^{2n}F^{2}(e_{a}, e_{b})=2ns\,.
\]
  \el
\pr
(1) The first expression follows since both $\mc{H}$ and $\mc{V}$ are $\fr{hol}(\nabla)$-invariant, hence $(\nabla_{X}Y_{\mc{H}})_{\mc{V}}=(\nabla_{X}Y_{\mc{V}})_{\mc{H}}$, see also \cite[Lemma 2.1.1]{ADS21} for a similar case.  For the second one we use Proposition \ref{2-valuedT} which implies that $T(X_{\mc{H}}, Y_{\mc{H}})_{\mc{V}}=T(X_{\mc{H}}, Y_{\mc{H}})=2F(X_{\mc{H}}, Y_{\mc{H}})\bar\xi\in\mc{V}$ and, moreover,  
\[
T(X_{\mc{H}}, Y_{\mc{V}})= \sum_{j=1}^s\eta_{j}(Y_{\mc{V}})T(X_{\mc{H}}, \xi_j)=2\sum_{j=1}^{s}\eta_{j}(Y_{\mc{V}})\phi(X_{\mc{H}})\,,
\]
for some $Y_{\mc{V}}=\sum_{j=1}^{s}\eta_{j}(Y_{\mc{V}})\xi_j\in\mc{V}$.  
The final relations in (1) are  now  immediate consequences of this expression.  For instance for $X, Y\in\mc{H}$ we have $Y_{\mc{V}}=0$, thus for any $Z\in\mc{V}$ we get
\[
g(\mc{A}_{X}Y, Z)=-F(X, Y)\sum_{j=1}^{s}\eta_{j}(Z)=-F(X, Y)\sum_{j=1}^{s}g(Z, \xi_j)=g(-F(X, Y)\bar{\xi}, Z)
\]
thus $\mc{A}_{X}Y=-F(X, Y)\bar{\xi}=-\frac{1}{2}T(X, Y)$.
Similarly,  
\[
\mc{A}_{X}\xi_j=(\nabla^{g}_{X}\xi_j)_{\mc{H}}=-(\phi X)_{\mc{H}}=-\phi X\,.
\]
(2) This is obtained using the same technique as in \cite[Lemma 2.2.1]{ADS21}, combining  Proposition \ref{2-valuedT} and Corollary \ref{XYvert}.  The established formula  is not used in the present work, but is included for completeness.   \\
(3) Recall that in terms of a $\phi$-basis we may assume that  $\mc{H}={\rm span}\{E_k, \phi(E_k) : 1\leq k\leq n\}$ and $\mc{V}={\rm span}\{\xi_\ell : 1\leq\ell\leq s\}$. Thus by part (1),  
we have $\mc{A}_{X}e_a=-F(X, e_a)\bar{\xi}$ for all $a\in\{1, \ldots, 2n\}$ and $X\in\mc{H}$. 
Moreover, since $F(X, E_i)=g(X, \phi(E_i))$ and $F(X, \phi(E_i))=-g(X, E_i)$ for all $i\in\{1, \ldots, n\}$,   
we obtain
\begin{eqnarray*}
\sum_{a=1}^{2n}g(\mc{A}_{X}e_{a}, \mc{A}_{Y}e_{a})&=&\sum_{a=1}^{2n}F(X, e_a)F(Y, e_a)g(\bar{\xi}, \bar{\xi})=s\sum_{a=1}^{2n}F(X, e_a)F(Y, e_a)\\
&=&s\sum_{i=1}^{n}\Big\{F(X, E_i)F(Y, E_i)+F(X, \phi(E_i))F(Y, \phi(E_i))\Big\}\\
&=&s\sum_{i=1}^{n}\Big\{g(X, \phi(E_i))g(Y, \phi(E_i))+g(X, E_{i})g(Y, E_i)\Big\}=sg(X, Y)\,.
\end{eqnarray*}
The second relation is now direct; we have   $\mc{A}_{e_{a}}e_{b}=-F(e_{a}, e_{b})\bar{\xi}$,  for all $a, b\in\{1, \ldots, 2n\}$, and this gives 
\[
\sum_{a, b=1}^{2n}g(\mc{A}_{e_{a}}e_{b}, \mc{A}_{e_{a}}e_{b})=\sum_{a, b=1}^{2n}F^2(e_a, e_b)g(\bar{\xi}, \bar{\xi})=s\sum_{a, b=1}^{2n}
g(e_a, \phi(e_b))g(e_a, \phi(e_b))=2ns\,.
\]
\pro
\bc\label{corolCurv} Let $\pi : (M^{2n+s}, g)\to (N^{2n}, g')$ be the canonical submersion of an $\mc{S}$-manifold. The Riemannian Ricci tensors $\Ric^g$ and $\Ric^{g'}$ are related by
\begin{equation}\label{Riccig'}
\Ric^{g}(X, Y)=\pi^*\Ric^{g'}(X, Y)-2\,s\,g(X, Y)\,,\quad X, Y\in\mc{H}\,.
\end{equation}
Therefore, the Riemannian scalar curvatures $\Scal^g$, $\Scal^{g'}$ and  the scalar curvature $\Scal^{\nabla}$ of the characteristic connection $\nabla$ on $M$ are related by
\begin{equation}\label{Scalg'}
\Scal^{g'}=\Scal^{g} +2\,n\,s=\Scal^{\nabla}+8\,n\,s\,.
\end{equation}
\ec
\pr
The Riemannian Ricci tensors of the total space and the base space  are such that (see \cite[pp.~243-244]{Bes})  
\[
\Ric^{g}(X, Y)=\pi^*\Ric^{g'}(X, Y)-2g(\mc{A}_X, \mc{A}_Y)\,,
\]
where $g(\mc{A}_X, \mc{A}_Y)$ is defined as in Lemma \ref{previous}. Thus, equation (\ref{Riccig'})  follows immediately from the result in  the same lemma.  
Let $\{e_a\}$ be a $g$-orthonormal $\phi$-adapted frame of $M^{2n+s}$. 
Based on (\ref{Riccig'}), we obtain 
\[
\sum_{a=1}^{2n}\Ric^{g}(e_a, e_a)=\sum_{a=1}^{2n}\Ric^{g'}(\pi_{*}e_a, \pi_{*}e_a)-2\,s\sum_{a=1}^{2n}g(e_a, e_a)=\Scal^{g'}-4\,n\, s\,.
\]
However, we know that  (see the proof in Corollary \ref{scalar_curvature})
\[
\Scal^g=\sum_{a=1}^{2n}\Ric^{g}(e_a, e_a)+\sum_{\ell=1}^{s}\Ric^{g}(\xi_{\ell}, \xi_{\ell})=\sum_{a=1}^{2n}\Ric^{g}(e_a, e_a)+2\,n\,s\,,
\]
thus the given relation between the Riemannian  scalar curvatures  follows. The second equation in (\ref{Scalg'}) follows from the relation $\Scal^{\nabla}=\Scal^{g}-6\,n\,s$, proved in Corollary \ref{scalar_curvature}. 
\pro

Recall that we can relate $\Ric^g$ with the  transverse Ricci tensor  $\Ric^{\check{g}}$ by  (\ref{Rictrans2}). Since the transverse metric $\check{g}$ satisfies 
 $g|_{\mc{D}\times\mc{D}}=\check{g}$, see (\ref{TKahlerg}), a comparison of this identity with the relation  (\ref{Riccig'})  obtained above, gives that

\bc
\label{Rel_transverse}
Given an $\mc{S}$-manifold $(M^{2n+s}, \phi, \xi_i, \eta_j, g)$ and the canonical submersion $\pi : M^{2n+s}\to N^{2n}$,  the pull-back $\pi^*\Ric^{g'}(X, Y)$ for $X, Y\in\mc{H}=\mc{D}$ coincides with the transverse Ricci tensor corresponding  to the transverse K\"ahler metric $\check{g}$ on the  characteristic foliation.
\ec

\br {\textsf{(The Ricci tensor on the fiber)}}
Let $\Pi:=\pi^{-1}(x')$ be the fiber over an arbitrary point $x'\in N$, where $\pi : M^{2n+s}\to N^{2n}$  is the canonical submersion of an $\mc{S}$-manifold.
We will denote by $\hat{g}=g|_{\Pi}$ the metric obtained via $g$ by restriction. Since each $\Pi$ is totally geodesic,  the Ricci tensor $\Ric^g$ and $\Ric^{\hat{g}}$ are related by
\[
\Ric^{g}(U, V)=\Ric^{\hat{g}}(U, V)+g(\mc{A}U, \mc{A}V)\,,\quad g(\mc{A}U, \mc{A}V):=\sum_{a=1}^{2n}g(\mc{A}_{e_a}U, \mc{A}_{e_{a}}V)
\]
for all vertical vector fields $U, V\in\mc{V}$ (see \cite{Bes}). Let $U=\sum_{i=1}^{s}\eta_{i}(U)\xi_i\in\mc{V}$ be an arbitrary vertical vector field. 
By Lemma \ref{previous} we get $\mc{A}_{e_{a}}U=-\sum_{i=1}^{s}\eta_{i}(U)\phi(e_a)$ for all $a\in\{1, \ldots, 2n\}$.  Thus 
\[
\sum_{a=1}^{2n}g(\mc{A}_{e_a}U, \mc{A}_{e_{a}}V)=\sum_{i, j=1}^{s}\eta_{i}(U)\eta_{j}(V)\sum_{a=1}^{2n}g(\phi(e_a), \phi(e_a))=2n\sum_{i, j=1}^{s}\eta_{i}(U)\eta_{j}(V)\,. 
\]
This gives
\begin{equation}\label{fiberRic}
\Ric^{g}(U, V)=\Ric^{\hat{g}}(U, V)+2\,n\, \sum_{i, j=1}^{s}\eta_{i}(U)\eta_{j}(V)\,.
\end{equation}
\er

\bex
Let $(N^{4}, g', J)$ be a simply connected compact K\"ahler-Einstein manifold with $\Scal^{g'}=32$,  and  assume that $N$  is Hodge, i.e., $\Omega\in H^{2}(N;\Z)$. 
By \cite[Theorem 2.5]{Blair03} there exists a circle bundle $\pi : M^{5}\to N^{4}$ together with a connection 1-form $\eta$ on $M^5$ such that $\dd\eta=\pi^*\Omega$. Moreover, $M$ should be compact, hence by Theorem \ref{BS70} the total space has a  regular Sasakian structure   $(\phi, \xi, \eta, g)$.
Therefore, we may identify $\pi$ with the $\BLY$-fibration for $s=1$, i.e., $\pi=\pi_1$.  By Example \ref{exampleSak} we know that   $\Scal^g=28$. Because  $n=2$ and $s=1$ we  can thus  confirm the relation $\Scal^{g'}=\Scal^{g} +2\,n\,s$, presented above.
\eex

\bex\label{U2again}
Let us consider the Lie group $\U(2)$ endowed with the left-invariant $\mc{S}$-structure $(\phi, \xi_i, \eta_j, g)$ mentioned in Example \ref{U2exam}.
Recall that  $n=1$ and  $s=2$.  The $\mc{S}$-structure is regular and the $\BLY$-submersion is given by 
\[
\Tg^2\longrightarrow \U(2)\overset{\pi_2}{\longrightarrow}  \CP^1\cong\Ss^2=\U(2)/\Tg^2\,,
\]
where $\Tg^2=\U(1)\times\U(1)$ is the maximal torus in $\U(2)$ (see also \cite{TK07}).
Hence the base space is the  K\"ahler manifold  $\CP^1$, which is  an isotropy irreducible Hermitian symmetric space
 endowed with its unique $\U(2)$-invariant Einstein metric which is K\"ahler (Fubini-Study metric).  
The submersion metric $g'$ induced by $\pi_2$ is $\U(2)$-invariant, hence by Schur's lemma it  coincides with  the Fubini-Study metric (up to scale).
The relation  $\Scal^{g'}=\Scal^{g} +2\,n\,s$ determines the normalization of the induced metric;  by  Example \ref{U2exam} we have $\Scal^g=4$
 and since $n=1$, $s=2$, we need $\Scal^{g'}=8$. 
Indeed,  for $\CP^n$, by \cite[9.81]{Bes} we know that the Fubini-Study metric  has Einstein constant  $2(n+1)$, thus the scalar curvature  of $\CP^n$ equals $4n(n+1)$,
so for $\CP^1$ we get scalar curvature $8$. Therefore, $g'$ is exactly the Fubini-Study metric (so we can remove the ambiguity ``up to scale''). 
Regarding the Ricci tensor on the fiber, the torus $\Tg^2$ is flat and  hence we have $\Ric^{\hat{g}}=0$, where $\hat{g}$ is the restriction of $g$ on $\Tg^2$.  Thus one may use    (\ref{fiberRic}) to confirm the relations
$\Ric^{g}(\xi_i, \xi_j)=2$ for $i, j\in\{1, 2\}$, presented in Example \ref{U2exam}.
\eex

   \bt\label{KEcor}
   Let  $(M^{2n+s}, \phi, \xi_i, \eta_j, g)$  be an $\mc{S}$-manifold and let  $\pi : M^{2n+s}\to N^{2n}$ be the canonical  $\fr{u}(n)$-submersion over the K\"ahler manifold $(N^{2n}, g', \Ja)$.  
   Suppose that the base space $N^{2n}$ is K\"ahler-Einstein with Einstein constant $\lambda$. Then the total space $M^{2n+s}$ is $\eta$-Einstein with $\al, \beta$ given by
   $\al=\lambda-2s$ and $\beta=2n-(\lambda-2s)$, respectively. Conversely, suppose that $(M^{2n+s}, \phi, \xi_i, \eta_j, g)$ is  an $\eta$-Einstein $\mc{S}$-manifold.  Then, the base space $(N^{2n}, g', \Ja)$ of   the canonical $\fr{u}(n)$-submersion $\pi$   is K\"ahler-Einstein with Einstein constant $\lambda=\al+2\, s$ (and  hence it is   Ricci-flat  if and only if $\al+2\,s=0$).
   \et
   \pr
 The result is essentially a consequence of Theorem \ref{Main_Corres},   since by Corollary \ref{Rel_transverse} one deduces that 
   the ``K\"ahler-Einstein condition'' $\Ric^{g'}=\lambda\cdot g'$  can be equivalently  replaced   by the condition that the $\mc{S}$-manifold $M^{2n+s}$ is
 transverse K\"ahler-Einstein with Einstein constant $\lambda$.  For completeness, let us use Corollary \ref{corolCurv} to describe the details.  The metrics $g$ on $M$ and $g'$ on $N$ are related by 
 \[
 g=\pi^*g'+\sum_{j=1}^{s}\eta_{j}\otimes\eta_j
 \]
  and by assumption we have $\Ric^{g'}=\lambda\cdot g'$. Thus, by   (\ref{Riccig'})   it follows that $\Ric^{g}(X, Y)=(\lambda-2s)g(X, Y)$, 
 for all $X, Y\in\mc{D}$.
 Similarly with the proof of    Theorem \ref{Main_Corres}, since we also have $\Ric^{g}(\xi_k, \xi_\ell)=2n$ for all $k, \ell\in\{1, \ldots, s\} $, 
  it follows that we may expresses $\Ric^g$ as 
 \[
 \Ric^{g}=(\lambda-2s)g+\beta\sum_{i}\eta_{i}\otimes\eta_{i}+2n\sum_{i\neq j}\eta_{i}\otimes\eta_{j}\,,
 \]
with $\beta=2n-(\lambda-2s)$. \\
\noindent Conversely, if $(M^{2n+s}, \phi, \xi_i, \eta_j, g)$  is an $\eta$-Einstein $\mc{S}$-manifold, then for the canonical submersion $\pi$  by the relations  (\ref{Ricgal}) and (\ref{Riccig'}) it follows that
      \[
      \pi^*\Ric^{g'}(X, Y)=(\al+2\,s)g(X, Y)=(\al+2\,s)\pi^*g'(X, Y)
      \]
      for all $X, Y\in\mc{H}=\mc{D}$, therefore  the base space is K\"ahler-Einstein with Einstein constant $\lambda=\al+2\, s$. 
   \pro

\bex
Consider the Lie group $\U(2)$ and the $\BLY$-fibration $\pi_2 : \U(2)\to\CP^1$ discussed above. The base space is K\"ahler-Einstein with Einstein constant $\lambda=4$, hence in this case we get $(\lambda-2s)=0$. Therefore $\U(2)$ endowed with its left-invariant $\mc{S}$-structure is $\eta$-Einstein with $\al=(\lambda-2s)=0$ and $\beta=2\cdot 1=2$.
 Based on the $\eta$-Einstein condition we may now confirm the matrix presentation of the Riemannian Ricci tensor presented in Example \ref{U2exam}, i.e.,
\begin{eqnarray*}
\Ric^{g}(\xi_1, \xi_1)&=&\beta\sum_{i=1}^{2}\eta_{i}(\xi_1)\eta_{i}(\xi_1)=2\big(\eta_1(\xi_1)\eta_{1}(\xi_1)+\eta_{2}(\xi_1)\eta_{2}(\xi_1)\big)
=2(1+0)=2\,,\\
\Ric^{g}(\xi_2, \xi_2)&=&2\big(\eta_1(\xi_2)\eta_{1}(\xi_2)+\eta_{2}(\xi_2)\eta_{2}(\xi_2)\big)=2\,,\\
\Ric^{g}(\xi_1, \xi_2)&=&(\al+\beta)\sum_{1\leq i\neq j\leq 2}\eta_{i}(\xi_1)\eta_{j}(\xi_2)
=2\big(\eta_{1}(\xi_1)\eta_{2}(\xi_2)+\eta_{2}(\xi_1)\eta_1(\xi_2)\big)=2
\end{eqnarray*}
and $\Ric^{g}(\xi_2, \xi_1)=2$, as well. \\  
Conversely, we just proved that  $\U(2)$ is $\eta$-Einstein. Hence, according to the theorem above $\CP^1$ should be K\"ahler-Einstein with Einstein constant $(\al+2\,s)$, which for $\al=0$ and $s=2$ gives the desired $4$.
\eex

%

\end{document}